\documentstyle[amssymb]{article}
\begin{document}
\noindent {\large \bf A Laplace operator and harmonics
on the quantum\\ complex vector space}
 \vskip 10 pt

{N. Z. Iorgov and A. U. Klimyk\footnote{Electronic mail:
aklimyk@bitp.kiev.ua}}

{\it Institute for Theoretical Physics, Kiev 03143, Ukraine}

\vskip  40 pt

\begin{abstract}
The aim of this paper is to study the $q$-Laplace operator and
$q$-harmonic polynomials on the quantum complex vector space
generated elements $z_i,w_i$, $i=1,2,\cdots ,n$, on which
the quantum group $GL_q(n)$ (or $U_q(n)$) acts. The
$q$-harmonic polynomials are defined as solutions of the equation
$\Delta _qp=0$, where $p$ is a polynomial in $z_i,w_i$,
$i=1,2,\cdots ,n$, and the $q$-Laplace operator $\Delta _q$ is
determined in terms of $q$-derivatives. The $q$-Laplace operator
$\Delta _q$ commutes with the action of $GL_q(n)$. The projector
${\sf H}_{m,m'}: {\cal A}_{m,m'}\to {\cal H}_{m,m'}$ is
constructed, where ${\cal A}_{m,m'}$ and ${\cal H}_{m,m'}$ are the
spaces of homogeneous (of degree $m$ in $z_i$ and of degree $m'$
in $w_i$) polynomials and homogeneous $q$-harmonic polynomials,
respectively. By using these projectors, a $q$-analogue of the
classical zonal spherical and associated spherical harmonics are
constructed. They constitute an orthogonal basis of ${\cal
H}_{m,m'}$. A $q$-analogue of separation of variables is given.
The quantum algebra $U_q({\rm gl}_n)$, acting on ${\cal
H}_{m,m'}$, determines an irreducible representation of $U_q({\rm
gl}_n)$. This action is explicitly constructed. The results of the
paper lead to the dual pair $(U_q({\rm sl}_2), U_q({\rm gl}_n))$
of quantum algebras.
\end{abstract}

\vskip 35 pt

\noindent {\sf I. INTRODUCTION}
\medskip

Laplace operators, harmonic polynomials, and related separations
of variables of the classical analysis are of a great importance
for mathematical and theoretical physics. They are closely related
to the rotation groups $SO(n)$ (if we deal with Euclidean space)
and to the unitary groups $U(n)$ (if we deal with the complex
vector space) (see, for example, Ref. 1, Chaps. 10 and 11). In this
paper we are interested in a $q$-analogue of Laplace operators,
harmonic
polynomials, and related separations of variables on complex
spaces.

Harmonic polynomials on the $n$-dimensional complex vector space
are defined by the equation $\Delta p=0$, where $\Delta$ is the
Laplace operator $\sum _{i=1}^n \partial^2/{\partial z_i \partial
{\bar z}_i}$ and $p$ belongs to the space ${\cal R}$ of
polynomials in $z_1,\cdots ,z_n, {\bar z}_1,\cdots ,{\bar z}_n$ on
the complex space ${\Bbb C}^n$. The space ${\cal H}$ of all
harmonic polynomials on ${\Bbb C}^n$ decomposes as a direct sum of
the subspaces ${\cal H}_{m,m'}$ of homogeneous harmonic
polynomials of degree $m$ in $z_1,\cdots ,z_n$ and of degree $m'$
in ${\bar z}_1,\cdots ,{\bar z}_n$: ${\cal H}=\bigoplus
_{m,m'=0}^\infty {\cal H}_{m,m'}$. The Laplace operator $\Delta$
on ${\Bbb C}^n$ commutes with the natural action of the unitary
group $U(n)$ on the space ${\Bbb C}^n$. This means that the
subspaces ${\cal H}_{m,m'}$ are invariant with respect to $U(n)$.
The irreducible representation $T_{m,m'}$ of the group $U(n)$ with
highest weight $(m,0,\cdots ,0,-m')$ is realized on ${\cal
H}_{m,m'}$.

The equation $\Delta p=0$ permits solutions in separated variables
on the space ${\cal H}_{m,m'}$. In other words, there exist
different coordinate systems (spherical, polyspherical, etc.) on
${\Bbb C}^n$ and for each of them it is possible to find the
corresponding basis of the space of solutions of the equation
$\Delta p=0$ consisting of products of functions depending on
separated variables (see Ref. 2 for the general theory of
separation of variables). To different coordinate systems there
correspond different separations of variables. From the other
side, to different coordinate systems there correspond different
chains of subgroups of the group $U(n)$ (see Ref. 1, Chap. 11, for
details of this correspondence). The bases of the space ${\cal
H}_{m,m'}$ in separated variables consist of products of Jacobi
polynomials multipled by $r^{m+m'}$ (different sets of Jacobi
polynomials for different separations of variables), where $r$ is
the radius. These polynomials (considered only on the unit sphere
$S^{n-1}_{\Bbb C}$ in ${\Bbb C}^n$) are matrix elements of the
class 1 (with respect to the subgroup $U(n-1)$) irreducible
representations $T_{m,m'}$ of $U(n)$ belonging to zero column (see
Ref. 1, Chap. 11).

Many new directions of mathematical physics are related to quantum
group and noncommutative geometry. It is natural to generalize
the above-described theory to noncommutative spaces. Such
generalizations can be of a great importance for further
development of some branches of mathematical and theoretical
physics related to noncommutative geometry.

The aim of this paper is to construct a $q$-deformation of the
above-described classical theory. In the $q$-deformed case,
instead of ${\Bbb C}^n$ we take the quantum complex vector space.
It is defined by the associative algebra ${\cal A}$ generated by
the elements $z_1,\cdots ,z_n,w_1,\cdots ,w_n$ satisfying a
certain natural defining relations. The elements $z_1,\cdots ,z_n$
play a role of Cartesian coordinates of ${\Bbb C}^n$ and
$w_1,\cdots ,w_n$ play a role of ${\bar z}_1,\cdots ,{\bar z}_n$.

The $q$-Laplace operator $\Delta _q$ on ${\cal A}$ is defined in
terms of $q$-derivatives (see formula (17) below). The quantum
group $U_q(n)$ play a role of the unitary group $U(n)$ in the
$q$-deformed case. It will be convenient for us to use the quantum
algebra (that is, the quantized universal enveloping algebra)
$U_q({\rm gl}_n)$ instead of the quantum group $U_q(n)$. The
$q$-harmonic polynomials on the quantum complex vector space are
defined as elements $p$ of the algebra ${\cal A}$ (that is,
polynomials in $z_1,\cdots ,z_n,w_1,\cdots ,w_n$) for which
$\Delta _q p=0$. By using the quantum algebra $U_q({\rm gl}_n)$ we
construct for $q$-harmonic polynomials a theory similar to the
theory for classical harmonic polynomials. We
construct the projector ${\sf H}_{m,m'}: {\cal A}_{m,m'}\to {\cal
H}_{m,m'}$, where ${\cal A}_{m,m'}$ and ${\cal H}_{m,m'}$ are the
subspaces of homogeneous (of degree $m$ in $z_1,\cdots ,z_n$ and
of degree $m'$ in $w_1,\cdots ,w_n$) polynomials in ${\cal A}$ and
in the space ${\cal H}$ of all $q$-harmonic polynomials from
${\cal A}$, respectively. Using these projectors we can make
different calculations in ${\cal H}_{m,m'}$. In this way, zonal
spherical and associated spherical polynomials can be calculated.
The associated spherical polynomials of ${\cal H}_{m,m'}$
constitute an orthogonal basis of this space. Here we obtain a
$q$-analogue of the spherical separation of coordinates. We show
that the natural action of the algebra $U_q({\rm gl}_n)$ on the
quantum complex vector space realizes on the space ${\cal
H}_{m,m'}$ the irreducible representation of this algebra with
highest weight $(m,0,\cdots ,0,-m')$. Note that restrictions of
zonal spherical and associated spherical polynomials from ${\cal
H}_{m,m'}$ to the quantum sphere in the quantum complex vector
space coincide with matrix elements of irreducible representations
$T_{m,m'}$ of the quantum group $U_q(n)$ corresponding to zero
column (the latter matrix elements were calculated in Ref. 3; see
also Ref. 4). Some our formulas coincide with formulas of Ref. 3.
However, no Laplace operator and no $q$-harmonic polynomials are
presented in Ref. 3.

Note that this paper is an extension of the results of our
previous paper (see Ref. 5) (where we studied $q$-Laplace operator
and $q$-harmonic polynomials on the quantum real vector space) to
the case of quantum complex vector space. It is well known that in
the classical case, the theory of Laplace operators and harmonic
polynomials on ${\Bbb C}^n$ can be reduced to the corresponding
theory for the real space ${\Bbb R}^{2n}$ (see Ref. 1, Chap. 11).
It is not the case for the quantum spaces. The reason is that the
quantum complex vector space cannot be obtained from the quantum
real vector space in the same way as in the classical case.

Everywhere below we suppose that $q$ is not a root of unity. Under
considering a scalar product on the spaces ${\cal A}$ and ${\cal
H}$ we assume that $q$ is a positive real number. By $[a]$, $a\in
{\Bbb C}$, we denote the so called $q$-number defined as
$$
[a]=\frac{q^a-q^{-a}}{q-q^{-1}}.
$$
\medskip

\noindent {\sf II. THE QUANTUM ALGEBRA $U_q({\rm gl}_n)$ AND THE
QUANTUM VECTOR SPACE}
\medskip

The Drinfeld--Jimbo quantum algebra $U_q({\rm gl}_n)$ is generated
by the elements $k^{1/2}_i\equiv q^{h_i/2}$, $k^{-1/2}_i\equiv
q^{-h_i/2}$, $i=1,2,\cdots ,n$, and $e_j$, $f_j$, $j=1,2,\cdots
,n-1$, satisfying the relations
$$
k_ik^{-1}_i=k^{-1}_ik_i=1,  \ \ \ k_ik_j=k_jk_i,\ \ \
k_ie_jk^{-1}_i=q^{a_{ij}}e_j, \ \ \ k_if_jk^{-1}_i=q^{-a_{ij}}f_j,
$$
$$
[e_i,f_j]\equiv e_if_j-f_je_i=\delta _{ij}
\frac{k_i k^{-1}_{i+1}-k^{-1}_i k_{i+1}}{q-q^{-1}},
$$      $$
e^2_ie_{i\pm 1}-(q+q^{-1})e_ie_{i\pm 1}e_i+e_{i\pm 1}e^2_i=0,
$$       $$
f^2_if_{i\pm 1}-(q+q^{-1})f_if_{i\pm 1}f_i+f_{i\pm 1}f^2_i=0,
$$      $$
[e_i,e_j]=[f_i,f_j]=0,\ \ \ \ \  |i-j|>1,
$$
+where $a_{ii}=1$, $a_{i,i-1}=a_{i-1,i}=-1$ and $a_{ij}=0$
otherwise (see, for example, Ref. 6, Chap. 6).

The algebra $U_q({\rm gl}_n)$ is a Hopf algebra, and the Hopf
algebra operations (comultiplication $\Delta$, counit
$\varepsilon$ and antipode $S$) are given by the formulas
$$
\Delta (k_i^{\pm 1})=k_i^{\pm 1}\otimes k_i^{\pm 1},\ \ \ \
\Delta (e_i)=e_i\otimes k^{-1/2}_ik^{1/2}_{i+1}+
k^{1/2}_ik^{-1/2}_{i+1}\otimes e_i,
$$     $$
\Delta (f_i)=f_i\otimes k^{-1/2}_ik^{1/2}_{i+1}+
k^{1/2}_ik^{-1/2}_{i+1}\otimes f_i,\ \ \ \
\varepsilon (k_i)=1,\ \ \ \
\varepsilon (e_i)=\varepsilon (f_i)=0,
$$      $$
S(k_i)=k_i^{-1},\ \ \ \
S(e_i)=-q^{-1}e_i,\ \ \ \ \
S(f_i)=- q f_i .
$$

The group $GL(n,{\Bbb C})$ and its Lie algebra ${\rm gl}(n,{\Bbb
C})$ act linearly on the $n$-dimensional complex vector space.
Similarly, the quantum group $GL_q(n,{\Bbb C})$ and the algebra
$U_q({\rm gl}_n)$ acts on the quantum (noncommutative) analogue of
the complex vector space. This quantum space is determined by the
algebra of polynomials ${\cal A}\equiv {\Bbb C}_q[z_1,\cdots
,z_n,$ $w_1, \cdots ,w_n]$ (see Ref. 7). This algebra is the
associative algebra generated by elements $z_1,z_2,\cdots ,z_n,$
$w_1,w_2, \cdots ,w_n$ satisfying the defining relations
$$
z_iz_j=qz_jz_i,\ \ \ \ w_iw_j=q^{-1}w_jw_i,\ \ \ \ i<j, \eqno (1)
$$   $$
w_jz_i=qz_iw_j,\ \ \ \ i\ne j,\ \ i,j=1,2,\cdots ,n,  \eqno (2)
$$   $$
w_kz_k=z_kw_k+(1-q^2)\sum _{s=1}^{k-1} z_sw_s .     \eqno (3)
$$
The elements $w_1,\cdots ,w_n$ play a role of ${\overline z_1},
\cdots ,{\overline z_n}$ in the classical analysis.

A $*$-operation can be defined on the algebra ${\cal A}$ turning
it into a $*$-algebra. This $*$-operation is uniquely determined
by the relations
$z_i^*=w_i$, $w_i^*=z_i$, $i=1,2,\cdots , n$.
The compact quantum group $U_q(n)$ acts on this $*$-algebra.

Note that the relations (3) are equivalent to the following ones:
$$
z_kw_k=w_kz_k-(1-q^2) \sum _{s=1}^{k-1} q^{2(k-s-1)}w_sz_s . \eqno
(4)
$$

The set of all monomials
$$
z_1^{r_1}z_2^{r_2}\cdots z_n^{r_n} w_1^{s_1}w_2^{s_2}\cdots
w_n^{s_n}, \ \ \ \ r_j,s_j=0,1,2,\cdots , \eqno(5)
$$
form a basis of the algebra ${\cal A}$ (see Ref. 8). The set
$$
w_1^{r_1}w_2^{r_2}\cdots w_n^{r_n} z_1^{s_1}z_2^{s_2}\cdots
z_n^{s_n}, \ \ \ \ r_j,s_j=0,1,2,\cdots , \eqno(6)
$$
also form a basis of this algebra.

The vector space of the algebra ${\cal A}$ can be represented as a
direct sum of the vector subspaces ${\cal A}_{m,m'}$
consisting of homogeneous polynomials of homogeneity degree $m$
in $z_1,z_2,\cdots ,z_n$ and of homogeneity degree $m'$ in
$w_1,w_2,\cdots ,w_n$, $m,m'=0,1,2,\cdots $:
$$
{\cal A}=\bigoplus _{m=0}^\infty \bigoplus _{m'=0}^\infty {\cal
A}_{m,m'} . \eqno (7)
$$
We have the linear space isomorphism
$$
{\cal A}\simeq {\cal A}_z\otimes {\cal A}_w,
$$
where the associative algebra ${\cal A}_z$ (the associative
algebra ${\cal A}_w$) is a subalgebra of ${\cal A}$ coinciding
with $\bigoplus _{m=0}^\infty {\cal A}_{m,0}$ (respectively, with
$\bigoplus _{m'=0}^\infty {\cal A}_{0,m'}$).

We can define an action of the algebra $U_q({\rm gl}_n)$ on the
vector space ${\cal A}$. To determine this action we give the
action of $U_q({\rm gl}_n)$ on $z_j$ and $w_j$ by the
formulas${}^8$
$$
k_i \triangleright z_j=q^{\delta _{ij}} z_j,\ \ \ \ e_i
\triangleright z_j=\delta _{j,i+1}z_{j-1}, \ \ \ \ f_i
\triangleright z_j=\delta _{j,i}z_{j+1} , \eqno(8)
$$      $$
k_i \triangleright w_j=q^{-\delta _{ij}} w_j,\ \ \ \ e_i
\triangleright w_j=-\delta _{j,i}q^{-1}w_{j+1}, \ \ \ \ f_i
\triangleright w_j=-\delta _{j,i+1}q w_{j-1} . \eqno(9)
$$
and extend it to ${\cal A}$ by using the comultiplication, that
is, by means of the relation
$$
X \triangleright (p_1 p_2)=\sum (X_{(1)}\triangleright p_1)
(X_{(2)}\triangleright p_2),
$$
where $\Delta (X)=\sum X_{(1)}\otimes X_{(2)}$ (in the Sweedler
notation), and linearity.

This action of the algebra $U_q({\rm gl}_n)$ on the vector space
${\cal A}$ determines a representation of $U_q({\rm gl}_n)$ on
this space (we denote it by $L$). Evidently, the subspaces ${\cal
A}_{m,m'}$ are invariant with respect to this action. Therefore,
$L$ determines representations of $U_q({\rm gl}_n)$ on these
subspaces, which are denoted by $L_{m,m'}$. We have $L=\bigoplus
_{m,m'=0}^\infty L_{m,m'}$.
\bigskip

\noindent {\sf 3. OPERATORS ON THE ALGEBRA ${\cal A}$}
\medskip

In order to introduce the $q$-Laplace operator on ${\cal A}$ and
to study $q$-harmonic polynomials we need some operators on the
linear space of the algebra ${\cal A}$. By $\gamma _i$ and
$\bar\gamma _i$ we denote the linear operators acting on monomials
as
$$
\gamma _i (z_1^{r_1}z_2^{r_2}\cdots z_n^{r_n}
w_1^{s_1}w_2^{s_2}\cdots w_n^{s_n})= q^{r_i}
z_1^{r_1}z_2^{r_2}\cdots z_n^{r_n} w_1^{s_1}w_2^{s_2}\cdots
w_n^{s_n},
$$ $$
\bar\gamma _i (w_1^{r_1}w_2^{r_2}\cdots w_n^{r_n}
z_1^{s_1}z_2^{s_2}\cdots z_n^{s_n})= q^{r_i}
w_1^{r_1}w_2^{r_2}\cdots w_n^{r_n} z_1^{s_1}z_2^{s_2}\cdots
z_n^{s_n}.
$$
Definition of the operators $\gamma^{-1} _i$ an $\bar\gamma^{-1}
_i$ is obvious.

By ${\hat z}_i$ and ${\breve z}_i$ we denote the linear operators
of multiplication by the element
$z_i$:
$$
{\hat z}_i (z_1^{r_1}z_2^{r_2}\cdots z_n^{r_n}
w_1^{s_1}w_2^{s_2}\cdots w_n^{s_n}) =z_i z_1^{r_1}z_2^{r_2}\cdots
z_n^{r_n} w_1^{s_1}w_2^{s_2}\cdots w_n^{s_n},
$$  $$
{\breve z}_i (z_1^{r_1}z_2^{r_2}\cdots z_n^{r_n}
w_1^{s_1}w_2^{s_2}\cdots w_n^{s_n}) = z_1^{r_1}z_2^{r_2}\cdots
z_n^{r_n} z_i w_1^{s_1}w_2^{s_2}\cdots w_n^{s_n}.
$$
The corresponding linear operators ${\hat w}_i$ and ${\breve w}_i$
are defined as
$$
{\hat w}_i (w_1^{r_1}w_2^{r_2}\cdots w_n^{r_n}
z_1^{s_1}z_2^{s_2}\cdots z_n^{s_n}) =w_i w_1^{r_1}w_2^{r_2}\cdots
w_n^{r_n} z_1^{s_1}z_2^{s_2}\cdots z_n^{s_n} ,
$$  $$
{\breve w}_i (w_1^{r_1}w_2^{r_2}\cdots w_n^{r_n}
z_1^{s_1}z_2^{s_2}\cdots z_n^{s_n}) = w_1^{r_1}w_2^{r_2}\cdots
w_n^{r_n} w_i z_1^{s_1}z_2^{s_2}\cdots z_n^{s_n} .
$$

We define on ${\cal A}$ the $q$-differentiations $\partial _i$ and
$\bar\partial _i$. The linear operators $\partial _i$ act as
$\partial _i p=0$ on monomials $p$ of the form (5) not containing
$z_i$ and as
$$
\partial _i ={\breve z}^{-1}_i \frac{\gamma _i-\gamma ^{-1}_i}{q-q^{-1}}
\eqno (10)
$$
on monomials containing $z_i$. The
$q$-differentiations $\bar\partial _i$ are linear operators acting
as $\bar\partial _i p=0$ on monomials $p$ of the form (6) not
containing $w_i$ and as
$$
\bar\partial _i ={\breve w}^{-1}_i \frac{\bar\gamma _i- \bar\gamma
^{-1}_i}{q-q^{-1}} \eqno (11)
$$
on monomials containing $w_i$.

The action formulas (8) and (9) mean that the multiplication
operators $\hat z_j$, $j=1,2,\cdots ,n$, and $\hat w_j$,
$j=1,2,\cdots ,n$, constitute tensor operators transforming under
the vector representation and under the contragredient to the
vector representation, respectively.

The actions (8) and (9) of $U_q({\rm gl}_n)$ on $z_j$ and $w_j$
determines its action on the operators $\partial_j$ and
$\bar\partial_j$:
$$
k_i \triangleright \partial_j=q^{-\delta _{ij}} \partial_j,\ \ \ \
e_i \triangleright \partial_j=-\delta _{j,i} q^{-1}
\partial_{j+1}, \ \ \ \ f_i \triangleright \partial_j=-\delta
_{j,i+1}q \partial_{j-1} , \eqno(12)
$$
$$
k_i \triangleright \bar\partial_j=q^{\delta _{ij}}
\bar\partial_j,\ \ \ \ e_i \triangleright \bar\partial_j= \delta
_{j,i+1}q^{-2}\bar\partial_{j-1}, \ \ \ \ f_i \triangleright
\bar\partial_j=\delta _{j,i}q^2 \bar\partial_{j+1} . \eqno(13)
$$
That is, the set $\bar\partial_j$, $j=1,2,\cdots,n$,
(respectively, the set $\partial_j$, $j=1,2,\cdots,n$) is a tensor
operator transforming under vector (respectively, contragredient
to vector) representation.

The operators $\partial _i$, $\bar\partial _i$, ${\hat z}_i$,
${\hat w}_i$ satisfy the relations, which will be presented by
means of the quantum $R$-matrix $R$ for the quantum algebra
$U_q({\rm gl}_n)$ (see, for example, Refs. 6, section 8.1, and 7
for definition of the $R$-matrix). Let ${\bf R}=P R$, where the
matrix $P$ permutes the spaces in the tensor product of two spaces
on which $R$-matrix acts. Then
\[
\mbox{\bf R}^{ij}_{kl}=q^{\delta_{ij}}\delta_{il}\delta_{jk}+
(q-q^{-1}) \delta_{ik}\delta_{jl} \theta(j-i),
\]
where $\theta(k)=1$ if $k>0$ and $\theta(k)=0$ if $k\le 0$. Its
inverse matrix is
\[
\bigl(\mbox{\bf
R}^{-1}\bigr)^{ij}_{kl}=q^{-\delta_{ij}}\delta_{il}\delta_{jk}-(q-q^{-1})
\delta_{ik}\delta_{jl} \theta(i-j).
\]
We also need the matrix $\Phi^{ij}_{kl}=\mbox{\bf R}^{ji}_{lk}
q^{2(i-l)}$, which satisfy the relations
\[
\sum_{j,l} \Phi^{ul}_{pj} \bigl(\mbox{\bf R}^{-1}\bigr)^{ji}_{lk}=
\sum_{j,l} \bigl(\mbox{\bf R}^{-1}\bigr)^{ul}_{pj}\Phi^{ji}_{lk}=
\delta_{up}\delta_{ik},
\]
\[
\sum_{l} \Phi^{li}_{lk}=\delta_{ik} q^{2(n-i)+1},\qquad
\sum_{k} \Phi^{jk}_{lk}=\delta_{jl} q^{2l-1}.
\]
The relations (1)--(3) rewritten for operators $\hat z_i$ and $\hat
w_i$ can be presented as
\[
\hat z_i \hat z_j =q^{-1} \mbox{\bf R}^{kl}_{ij} \hat z_k \hat z_l, \qquad
\hat w_i \hat w_j =q^{-1} \mbox{\bf R}^{ji}_{lk} \hat w_k \hat w_l, \qquad
\hat w_i \hat z_j =q \bigl(\mbox{\bf R}^{-1}\bigr)^{ik}_{jl} \hat z_k \hat
w_l,
\]
We also have the relations
\[
\partial_i \partial_j =q^{-1} \mbox{\bf R}^{ji}_{lk} \partial_k
\partial_l, \qquad \bar\partial_i \bar\partial_j =q^{-1}
\mbox{\bf R}^{kl}_{ij} \bar\partial_k \bar\partial_l, \qquad
\partial_i \bar\partial_j =q^{-1} \Phi^{ki}_{lj} \bar\partial_k \partial_l,
\]
\[
\partial_i \hat w_j=q \bigl(\mbox{\bf R}^{-1}\bigr)^{ji}_{lk}
\hat w_k \partial_l,\qquad
\bar\partial_i \hat z_j=q \Phi^{lk}_{ji} \hat z_k \bar\partial_l,
\]
\[
\partial_i \hat z_j= \gamma^{\mp 1}\delta_{ij}+({\bf R}^{\pm 1})^{ik}_{jl}
\hat z_k \partial_l,\qquad \bar\partial_i \hat w_j=
\bar\gamma^{\mp 1}\delta_{ij}+ ({\bf R}^{\pm 1})^{lj}_{ki} \hat
w_k \bar\partial_l,
\]
which can be represented in the form
\[
\partial_i \partial_j = q^{-1}\partial_j\partial_i, \qquad
\bar\partial_i \bar\partial_j = q\bar\partial_j\bar\partial_i,
\qquad i<j,
\]
\[
\bar\partial_i \partial_j = q\partial_j\bar\partial_i, \ \ i\neq
j,\qquad \bar\partial_i \partial_i = \partial_i \bar\partial_i+
(1-q^2)\sum_{k>i} \partial_k \bar\partial_k,
\]
\[
\partial_i \bar\partial_i = \bar\partial_i \partial_i+
(1-q^{-2})\sum_{k>i} q^{2(k-i)}\bar\partial_k \partial_k.
\]
\[
\partial_i \hat w_i=\hat w_i \partial_i,\qquad
\partial_i \hat w_j - q \hat w_j \partial_i=
(1-q^2) \hat w_i \partial_j,\qquad
\partial_j \hat w_i = q \hat w_i \partial_j, \qquad  i<j,
\]
\[
\bar\partial_i \hat z_i = \hat z_i \bar\partial_i,\qquad
\bar\partial_i \hat z_j - q^{-1} \hat z_j \bar\partial_i=
(1-q^{-2})q^{2(j-i)} \hat z_i \bar\partial_j,\qquad \bar\partial_j
\hat z_i = q^{-1} \hat z_i \bar\partial_j, \qquad  i<j,
\]
\[
\partial_i \hat z_j = \hat z_j \partial_i, \qquad
\bar\partial_i \hat w_j = \hat w_j \bar\partial_i, \qquad  i\neq j,
\]
\[
\partial_i \hat z_i=q \hat z_i \partial_i+ (q-q^{-1}) \sum_{k>i}
\hat z_k \partial_k+\gamma^{-1}= q^{-1} \hat z_i
\partial_i- (q-q^{-1}) \sum_{k<i} \hat z_k \partial_k+\gamma,
\]
\[
\bar\partial_i \hat w_i=q \hat w_i \bar\partial_i+ (q-q^{-1})
\sum_{k<i} \hat w_k \bar\partial_k+\bar\gamma^{-1}=q^{-1} \hat w_i
\bar\partial_i- (q-q^{-1}) \sum_{k>i} \hat w_k
\bar\partial_k+\bar\gamma,
\]
where $\gamma =\gamma_1 \gamma_2\cdots \gamma_n$ and $\bar\gamma
=\bar\gamma_1 \bar\gamma_2\cdots \bar\gamma_n$. From last two
lines, we obtain
\[
\sum_{k=1}^n \hat z_k
\partial_k=\{\gamma\}\equiv\frac{\gamma-\gamma^{-1}}{q-q^{-1}},
\qquad \sum_{k=1}^n \hat w_k \bar\partial_k=\{\bar\gamma\}\equiv
\frac{\bar\gamma-\bar\gamma^{-1}}{q-q^{-1}},
\]
We also have the relations
\[
\gamma {\hat z}_i=q {\hat z}_i \gamma, \qquad
\gamma {\hat w}_i= {\hat w}_i \gamma, \qquad
\bar \gamma {\hat z}_i= {\hat z}_i \bar\gamma, \qquad
\bar \gamma {\hat w}_i= q {\hat w}_i \bar\gamma,
\]
\[
\gamma {\partial}_i=q^{-1} {\partial}_i \gamma, \qquad
\gamma {\bar \partial}_i= {\bar \partial}_i \gamma, \qquad
\bar \gamma {\partial}_i= {\partial}_i \bar\gamma, \qquad
\bar \gamma \bar{\partial}_i= q^{-1} \bar{\partial}_i \bar\gamma.
\]
Note that
$$
\gamma p = q^m p,\qquad  \bar \gamma p=q^{m'} p,\qquad
p \in {\cal A}_{m,m'}. \eqno(14)
$$
To compare these relations with known from literature, we introduce
the operators $\partial'_i=\gamma\partial_i$,
$\bar\partial'_i=\bar\gamma^{-1}\bar\partial_i$. Then the
operators $\hat z_i$, $\hat w_i$, $\partial'_i$,
$\bar\partial'_i$, $i=1,\cdots,n,$ satisfy the relations from Ref.
9 which are known to be covariant with respect to $U_q({\rm
gl}_n)$.

Note that the above elements ${\hat z}_1, \cdots {\hat z}_n,
\partial' _1,\cdots ,\partial' _n$ generate the $q$-Weyl algebra,
that is, they satisfy the relations
$$
{\hat z}_i{\hat z}_j=q {\hat z}_j {\hat z}_i,\ \ \
\partial'_i \partial' _j=q^{-1}\partial' _j \partial'_i\ \  \ \ i<j,
\qquad
\partial' _i{\hat z}_j=q{\hat z}_j \partial' _i, \ \ \ \ \ i\ne j,
$$    $$
\partial' _i {\hat z}_i-q^2{\hat z}_i \partial'_i=
1+(q^2-1)\sum _{j>i}{\hat z}_j \partial'_j
$$
(the definition of the $q$-Weyl algebra see, for example, in Ref.
6, Chap. 12). Similarly, the elements ${\hat w}_1, \cdots {\hat
w}_n, \bar\partial' _1,\cdots ,\bar\partial' _n$ generate the
$q^{-1}$-Weyl algebra.

The operators
$D:=\sum _{k=1}^n {\hat z}_k\partial _k$ and ${\bar D}:=\sum
_{k=1}^n {\hat w}_k\bar\partial _k$
are called the {\it $q$-Euler operators}. The formula (7) gives
the decomposition of ${\cal A}$ into a direct sum of eigenspaces
of the operators $D$ and $\bar D$.

Let us show that the above relations for the operators
$\partial_i, \bar\partial _i, {\hat z}_i, {\hat w}_i$ determine
uniquely the formulas (10) and (11) for $\partial_i, \bar\partial
_i$. We use the action formulas $\partial_i 1=\bar\partial_i 1=0$,
take into account that $\hat z_i$, $\hat w_i$ act as the operators
of left multiplication on the basis elements (5) and (6),
respectively, and $\gamma$, $\bar\gamma$ are gradation operators
on ${\cal A}$ (see (14)). By means of commutation relations
between $\partial_i$ and $\hat w_j$, it is easy to obtain that
$\partial_i w_1^{s_1}w_2^{s_2}\cdots w_n^{s_n}=0$. To calculate
$\partial_i (z_i^{r_i} w_1^{s_1}w_2^{s_2}\cdots w_n^{s_n})$ with
$r_i>0$, we use the relation
\[
\partial _i \hat z_i= q \hat z_i \partial _i +
(q-q^{-1})\sum_{k>i} \hat z_k \partial_k+\gamma^{-1}.
\]
It gives $\partial_i (z_i^{r_i} w_1^{s_1}w_2^{s_2}\cdots
w_n^{s_n})=$ $[r_i] z_i^{r_i-1} w_1^{s_1}w_2^{s_2}\cdots
w_n^{s_n}$. Finally, we have the action formula
\[
\partial_i (z_1^{r_1}\cdots z_i^{r_i}\cdots z_n^{r_n}
w_1^{s_1}\cdots w_n^{s_n})=
\]
\[
=q^{r_i(r_{i+1}+\cdots+r_{n})}\partial_i (z_1^{r_1}\cdots
z_{i-1}^{r_{i-1}} z_{i+1}^{r_{i+1}} \cdots z_n^{r_n} z_{i}^{r_{i}}
w_1^{s_1}\cdots w_n^{s_n})=
\]
\[=
q^{r_i(r_{i+1}+\cdots+r_{n})} \hat z_1^{r_1}\cdots \hat
z_{i-1}^{r_{i-1}} \hat z_{i+1}^{r_{i+1}} \cdots \hat z_n^{r_n}
\partial_i (z_{i}^{r_{i}} w_1^{s_1}\cdots w_n^{s_n})=
\]
\[=
q^{r_{i+1}+\cdots+r_{n}}[r_i] z_1^{r_1}\cdots z_i^{r_i-1}\cdots
z_n^{r_n} w_1^{s_1}\cdots w_n^{s_n},
\]
which exactly coincides with the action (10). The action formula
for $\bar\partial_i$ is recovered in a similar way.

The action of the algebra $U_q({\rm gl}_n)$ on ${\cal A}\simeq
{\cal A}_z\otimes {\cal A}_w$, defined by formulas (8) and
(9), can be determined in terms of the operators $\partial _i$
and $\bar\partial _j$. We first note that the action of $U_q({\rm
gl}_n)$ on ${\cal A}_z$ is given by the operators
\[
L(k_i)=\gamma _i,\ \ \ \
 L(e_i)=q^{-1/2}(\gamma _i\gamma _{i+1})^{1/2}{\breve z}_i
 \partial _{i+1},\ \ \ \
 L(f_i)=q^{1/2}(\gamma _i\gamma _{i+1})^{-1/2}{\breve z}_{i+1}
 \partial _{i}
\]
and on ${\cal A}_w$ by the operators
$$
L(k_i)= \bar\gamma ^{-1}_{i},\ \
 L(e_i)=-q^{-3/2}(\bar\gamma _i\bar\gamma _{i+1})^{1/2}
 {\breve w}_{i+1}  \bar\partial _{i},\ \
 L(f_i)=-q^{3/2}(\bar\gamma _i\bar\gamma _{i+1})^{-1/2}
 {\breve w}_i  \bar\partial _{i+1}.
$$
Taking into account the comultiplication on $U_q({\rm gl}_n)$ the
action of $U_q({\rm gl}_n)$ on the linear space ${\cal A}\simeq
{\cal A}_z\otimes {\cal A}_w$ can be written as $L(k_i)=\gamma _i
\otimes \bar\gamma ^{-1}_{i}$ and
\[
L(e_i)=q^{-1/2}(\gamma _i\gamma _{i+1})^{1/2}{\breve z}_i
 \partial _{i+1}\otimes
(\bar\gamma_i \bar\gamma^{-1} _{i+1})^{1/2}-q^{-3/2}
(\gamma_i \gamma^{-1} _{i+1})^{1/2}
\otimes (\bar\gamma _i\bar\gamma _{i+1})^{1/2}
 {\breve w}_{i+1}  \bar\partial _{i},
\]\[
L(f_i)=q^{1/2}(\gamma _i\gamma _{i+1})^{-1/2}{\breve z}_{i+1}
 \partial _{i}\otimes
(\bar\gamma_i \bar\gamma^{-1} _{i+1})^{1/2}-q^{3/2}
(\gamma_i \gamma^{-1} _{i+1})^{1/2}
\otimes (\bar\gamma _i\bar\gamma _{i+1})^{-1/2}
 {\breve w}_i  \bar\partial _{i+1}.
\]
\medskip

\noindent {\sf IV. SQUARED $q$-RADIUS AND $q$-LAPLACE OPERATOR}
\medskip

The element
$$
Q=\sum _{i=1}^n z_iw_i=\sum _{i=1}^n q^{2(n-i)}w_iz_i \in {\cal
A}_{1,1} \eqno (15)
$$
of the algebra ${\cal A}$ is called the {\it squared $q$-radius}
on the quantum complex vector space. It is an important element in
${\cal A}$. One can check by a direct computation that $Q$ {\it is
invariant with respect to the representation $L_{1,1}$} (and hence
with respect to the representation $L$), that is,
$L(k^{\pm 1}_i)Q=Q$, $L(e_j)Q=0$ and $L(f_j)Q=0$.
Similarly, the element $Q^k\in {\cal A}_{k,k}$ is invariant with
respect to the representation $L_{k,k}$.

The squared $q$-radius $Q$ belongs to the center of the algebra
${\cal A}$, that is,
$Qz_i=z_iQ$, $Qw_i=w_iQ$, $i=1,2,\cdots ,n$.
We shall also use the elements
$$
Q_j=\sum _{i=1}^j z_iw_i=\sum _{i=1}^j q^{2(j-i)}w_iz_i,
$$
which are squared $q$-radiuses for the subalgebras ${\Bbb
C}_q[z_1,w_1,\cdots , z_j,w_j]$. They satisfy the relations${}^8$
$$
Q_jQ_i=Q_iQ_j,\ \ \ \   z_iw_i=Q_i-Q_{i-1},\ \ \ \
w_iz_i=Q_i-q^2Q_{i-1} ,
$$    $$
z_iQ_j=q^{-2}Q_jz_i,\ \ \  w_iQ_j=q^2Q_jw_i\ \ \ \ {\rm for}\ \ \ \
i>j,
$$    $$
z_iQ_j=Q_jz_i,\ \ \  w_iQ_j=Q_jw_i\ \ \ \ {\rm for}\ \ \ \ i\le j.
$$

It can be checked${}^8$ by direct computation that
$$
z^k_iw^k_i=Q^k_i\left( Q_{i-1}/Q_i;q^{-2}\right) _k,\ \ \
w^k_iz^k_i=Q^k_i\left( q^2Q_{i-1}/Q_i;q^{2}\right) _k,
\eqno (16)
$$
where
$$
(a;q)_s=(1-a)(1-aq)\cdots (1-aq^{s-1}).
$$

We consider on ${\cal A}$ the operator
$$
\Delta _q =\partial_1 \bar\partial _1 +
\partial _2\bar\partial _2+ \cdots +
\partial _n\bar\partial _n =\sum _{i=1}^n q^{2(i-1)}
\bar\partial _i  \partial _i , \eqno (17)
$$
which is called the $q$-{\it Laplace operator} on the quantum
complex vector space. Since $\gamma \Delta_q=q^{-1} \Delta_q
\gamma$ and $\bar\gamma \Delta_q=q^{-1} \Delta_q \bar\gamma$, then
$\Delta_q : {\cal A}_{m,m'}\to {\cal A}_{m-1,m'-1}$.

To the element (15) there corresponds the operator ${\hat Q}$
on ${\cal A}$ defined as
$$
{\hat Q}={\hat z}_1{\hat w}_1+{\hat z}_2{\hat w}_2+\cdots
+{\hat z}_n{\hat w}_n .
$$

{\it Proposition 1:} {\it The operators $\Delta_q$ and ${\hat Q}$
satisfy the relations
$$
\Delta_q {\hat Q}^k- {\hat Q}^k\Delta_q =q^{n-1}{\hat
Q}^{k-1}[k]\{ q^{k+n-1}\gamma\bar\gamma \} ,  \eqno (18)
$$     $$
\Delta_q ( Q^k)= q^{n-1}Q^{k-1}[k][k+n-1], \eqno (19)
$$
where
$$
\{ a\} =\frac{a-a^{-1}}{q-q^{-1}}
$$
and $[r]\equiv \{q^r\}$ is a $q$-number.}

{\it Proof:} First we prove the relation $[\Delta_q, {\hat Q}]
=q^{n-1}\{q^{n}\gamma\bar\gamma\}$. Using relations of section III
we derive
\[
\Delta_q \hat Q=\sum_{k,l}\partial_k \bar\partial_k \hat z_l \hat
w_l= \sum_{k,l,i,j} \partial_k \bigl(q^{-1} \Phi^{ji}_{lk} \hat
z_i \bar\partial_j \bigr) \hat w_l
\]
\[
=\sum_{k,l,i,j} q^{-1} \Phi^{ji}_{lk} \bigl(\delta_{ik}
\gamma^{-1}+\sum_{r,s} \mbox{\bf R}^{kr}_{is}\hat z_r \partial_s
\bigr) \bar\partial_j \hat w_l
\]
\[
=\sum_{k,l,j} q^{-1} \Phi^{jk}_{lk} \gamma^{-1} \bar\partial_j
\hat w_l+ \sum_{k,l,i,j,r,s} q^{-1} \Phi^{ji}_{lk} \mbox{\bf
R}^{kr}_{is} \hat z_r \partial_s  \bigl(\delta_{jl} \bar\gamma+
\sum_{u,p} \bigl(\mbox{\bf R}^{-1}\bigr)^{ul}_{pj}\hat w_p
\bar\partial_u \bigr)
\]
\[
=\sum_{l} q^{2l-2}\gamma^{-1} \bar\partial_l \hat w_l+
\sum_{i,r,s} q^{2(n-i)} \mbox{\bf R}^{ir}_{is}  \hat z_r
\partial_s \bar\gamma+ \sum_{k,i,r,s} q^{-1} \mbox{\bf
R}^{kr}_{is} \hat z_r \partial_s  \hat w_i \bar\partial_k.
\]
The third summand is equal to
\[
\sum_{k,i,r,s} q^{-1} \mbox{\bf R}^{kr}_{is} \hat z_r \partial_s
\hat w_i \bar\partial_k= \sum_{k,i,r,s} q^{-1} \mbox{\bf
R}^{kr}_{is} \hat z_r \bigl(q \sum_{u,p} \bigl(\mbox{\bf
R}^{-1}\bigr)^{is}_{pu} \hat w_u \partial_p \bigr)
\bar\partial_k={\hat Q} \Delta_q.
\]
Using explicit expressions for matrix elements of {\bf R} and
$\mbox{\bf R}^{-1}$ we have
\[
\sum_{i} q^{2(n-i)} \mbox{\bf
R}^{ir}_{is}=q^{2n-1}\delta_{rs},\qquad \sum_{l} q^{2l-2}
\bigl(\mbox{\bf R}^{-1}\bigr)^{ul}_{pl}=q^{-1}\delta_{pu} ,
\]
\[
\sum_{l} q^{2l-2} \bar\partial_l \hat w_l= \sum_{l} q^{2l-2}
\bigl(\bar\gamma+\sum_{u,p} \bigl(\mbox{\bf
R}^{-1}\bigr)^{ul}_{pl}\hat w_p \bar\partial_u \bigr)
\]
\[
=q^{n-1}[n]\bar\gamma+q^{-1}\sum_p \hat w_p
\bar\partial_p=q^{n-1}[n]\bar\gamma+q^{-1}\{\bar\gamma\}=
q^{n-1}\{q^n\bar\gamma\}
\]
Thus, $[\Delta_q, {\hat Q}]=q^{n-1}\gamma^{-1}\{
q^n\bar\gamma\}+q^{2n-1}\{\gamma\}
\bar\gamma=q^{n-1}\{q^n\gamma\bar\gamma\}$.
Now, it is easy to obtain (18) by induction if to use the relation
$\{q^r \gamma\bar\gamma\} \hat Q= \hat Q \{q^{r+2}
\gamma\bar\gamma\}$ and the explicit expression for $\{a\}$.
Acting by both sides of (18) on $1$ we obtain (19).

{\it Proposition 2:} {\it The operators $\Delta_q$ and ${\hat Q}$
commute with the action of the algebra $U_q({\rm gl}_n)$ on ${\cal
A}$, that is, with all operators of the representation $L$ of
$U_q({\rm gl}_n)$.}

{\it Proof:} It follows from (12) and (13) that
$k_i\triangleright\Delta_q=\Delta_q$, $e_j
\triangleright\Delta_q=0$ and $f_j \triangleright\Delta_q=0$. Now
using the comultiplication for $k_i$, $e_j$ and $f_j$, we obtain
the proposition for the $q$-Laplace operator. For ${\hat Q}$ the
proposition is proved similarly.
\bigskip

\noindent {\sf V. $q$-HARMONIC POLYNOMIALS}
\medskip

A polynomial $p\in {\cal A}$ is called $q$-{\it harmonic} if
$\Delta_q p=0$. The linear subspace of ${\cal A}$ consisting of
all $q$-harmonic polynomials is denoted by ${\cal H}$. Let
$$
{\cal H}_{m,m'}={\cal A}_{m,m'}\cap {\cal H}.
$$

{\it Proposition 3:} {\it The space ${\cal A}_{m,m'}$ can be
represented as the direct sum}
$$
{\cal A}_{m,m'}={\cal H}_{m,m'}\oplus Q{\cal A}_{m-1,m'-1}. \eqno
(20)
$$

{\it Proof:} First we prove that ${\cal H}_{m,m'}\cap Q{\cal
A}_{m-1,m'-1}=\{ 0\}$. If it is not true, then there exists
nonzero element $p\in {\cal H}_{m,m'}\cap Q{\cal A}_{m-1,m'-1}$.
Let $k$ be a maximal integer such that $p=Q^kp'$ with some nonzero
polynomial $p'$. Then it follows from $\Delta_q (p)=0$ and (18)
that
$$
0=\Delta_q (Q^kp')=Q^k \Delta_q
(p')+Q^{k-1}q^{n-1}[k][k+n-1+m+m'-2k] p'.
$$
Since $q^{n-1}[k][k+n-1+m+m'-2k]\ne 0$, then $p'$ can be divided
by $Q$. This is a contradiction. Thus, ${\cal H}_{m,m'}\cap Q{\cal
A}_{m-1,m'-1}=\{ 0\}$. Using this fact and the equality $\ker
\Delta_q={\cal H}_{m,m'}$, where $\Delta_q$ is considered only on
${\cal A}_{m,m'}$, we obtain the chain of inequalities
\[
\dim {\cal A}_{m,m'} - \dim \ker \Delta_q\ge \dim Q {\cal
A}_{m-1,m'-1}= \dim {\cal A}_{m-1,m'-1}\ge \dim\mathop{{\rm im}}
\Delta_q.
\]
The last inequality follows from the fact that $\Delta_q : {\cal
A}_{m,m'}\to {\cal A}_{m-1,m'-1}$. Now we take into account the
relation $\dim \ker \Delta_q + \dim\mathop{{\rm im}} \Delta_q=\dim
{\cal A}_{m,m'}$. Thus, in fact, the above inequalities are exact
equalities, and ${\cal A}_{m,m'}={\cal H}_{m,m'}\oplus Q{\cal
A}_{m-1,m'-1}$. Proposition is proved.

{\it Remark:} If $n=1$, then ${\cal A}$ consists of all
polynomials in commuting elements $z_1$ and $w_1$. In this case,
the space ${\cal H}$ of $q$-harmonic polynomials has a basis
consisting of the polynomials
$$
1, \ \ z_1^k,\ \ w_1^k,\ \ \ \ k=1,2,\cdots . \eqno (21)
$$

The decomposition (20) has also the following consequences:

{\it Corollary 1:} {\it If $p\in {\cal H}_{m,m'}$, then $p$ cannot
be represented as $p=Q^kp'$, $k\ne 0$, with some polynomial $p'$.}

{\it Corollary 2:} {\it The space ${\cal A}_{m,m'}$ decomposes
into the direct sum}
$$
{\cal A}_{m,m'}=\bigoplus _{j=0}^{\min(m,m')} Q^j{\cal
H}_{m-j,m'-j} .\eqno (22)
$$

{\it Corollary 3:} {\it For dimension of the space of $q$-harmonic
polynomials ${\cal H}_{m,m'}$ we have the formula}
$$
{\rm dim}\ {\cal H}_{m,m'}={(m+n-2)!(m'+n-2)!(m+m'+n-1)\over
(n-1)!(n-2)!m!m'!} .
$$

{\it Corollary 4:} {\it The space of $q$-harmonic polynomials
${\cal H}$ can be represented in the form of a direct sum}
$$
{\cal H}=\bigoplus _{m=0}^\infty \bigoplus _{m'=0}^\infty {\cal
H}_{m,m'} .
$$

Corollary 1 is a direct consequence of formula (20). Corollary 2
easily follows from repeated application of (20). Corollary 3 is
proved in the same way as in the classical case (see, for example,
Ref. 1, Chap. 10). For this we note that
$$
{\dim}\ {\cal A}_{m,m'}={(n+m-1)!(n+m'-1)!\over (n-1)!^2m!m'!} .
$$
Hence, for ${\dim}\ {\cal H}_{m,m'}={\dim}\ {\cal A}_{m,m'}-
{\dim}\ {\cal A}_{m-1,m'-1}$ we obtain the expression stated in
the corollary.
In order to prove Corollary 4 we note that
$$
{\cal A}=\bigoplus _{m\ge 0}\bigoplus _{m'\ge 0} \bigoplus
_{j=0}^p Q^j{\cal H}_{m-j,m'-j}= \bigoplus _{m\ge 0}\bigoplus
_{m'\ge 0}\left( {\cal H}_{m,m'}\bigoplus \Bigl( \bigoplus
_{j=1}^p Q^j{\cal H}_{m-j,m'-j}\Bigr) \right) ,
$$
where $p=\min(m,m')$. Now Corollary 4 follows from here and
Corollary 1.

{\bf Theorem 1:} {\it The linear space isomorphism ${\cal A}\simeq
{\Bbb C}[Q]\otimes {\cal H}$ is true, where ${\Bbb C}[Q]$ is the
space of all polynomials in $Q$.}

This theorem follows from Corollary 2.

The decomposition ${\cal A}\simeq {\Bbb C}[Q]\otimes {\cal H}$ is
a $q$-analogue of the theorem on separation of variables for Lie
groups in an abstract form${}^{10}$. It follows from this
decomposition that
$$
{\cal A}\simeq {\Bbb C}[Q]\otimes {\cal H}\simeq {\Bbb
C}[Q]\otimes \bigoplus _{m\ge 0}\bigoplus _{m'\ge 0} {\cal
H}_{m,m'} = \bigoplus _{m\ge 0}\bigoplus _{m'\ge 0} \left( {\Bbb
C}[Q]\otimes {\cal H}_{m,m'} \right) . \eqno (23)
$$

Since the subspaces ${\cal A}_{m,m'}$ are invariant with respect
to the action of the algebra $U_q({\rm gl}_n)$, it follows from
Proposition 2 for $\Delta_q$ that the subspace ${\cal H}_{m,m'}$
is invariant with respect to the representation $L_{m,m'}$ of
$U_q({\rm gl}_n)$. We denote the restriction of this
representation to ${\cal H}_{m,m'}$ by $T_{m,m'}$. It follows from
Proposition~2 for $Q$ and from (22) that
$$
L_{m,m'}=\bigoplus _{j=0}^{\min (m,m')} T_{m-j,m'-j} . \eqno (24)
$$

{\it Proposition 4:} {\it The representations $T_{m-j,m'-j}$ of
$U_q({\rm gl}_n)$ in (24) are irreducible with highest weights
$(m-j,0,\cdots ,0,-m'+j)$, respectively.}

{\it Proof:} Let us show that the representation $L_{m,0}=T_{m,0}$
in the space of holomorphic polynomials ${\cal A}_{m,0}$ is
irreducible with highest weight $(m,0,\cdots ,0)$. In fact, a
direct calculation shows that the monomials $z_1^{m_1}\cdots
z_n^{m_n}$, $m_1+\cdots +m_n=m$, are weight vectors of this
representation. The highest weight vector coincides with $z_1^m$.
Therefore, the irreducible representation with highest weight
$(m,0,\cdots ,0)$ is a subrepresentation of $L_{m,0}=T_{m,0}$.
Since their dimensions coincide, $L_{m,0}=T_{m,0}$ is an
irreducible representation with highest weight $(m,0,\cdots ,0)$.
It can be proved in the same way that the representation
$L_{0,m'}=T_{0,m'}$ in the space of polynomials ${\cal A}_{0,m'}$
is irreducible with highest weight $(0,\cdots ,0, -m')$.

Now we can prove the proposition by the induction. Assume that the
proposition is true for the representations $T_{m-1-j,m'-1-j}$
which are contained in the decomposition
$$
L_{m-1,m'-1}=\bigoplus _{j=0}^{\min (m-1,m'-1)} T_{m-1-j,m'-1-j} .
\eqno (25)
$$
Note that since
${\cal A}_{m,m'}={\cal H}_{m,m'}\oplus Q{\cal A}_{m-1,m'-1}$, then
$L_{m-1,m'-1}$ is a subrepresentation in $L_{m,m'}$ and
$$
\dim {\cal A}_{m-1,m'-1}=\dim L_{m-1,m'-1} =\sum _{j=0}^{\min
(m-1,m'-1)} \dim T_{m-1-j,m'-1-j}.
$$
The space ${\cal A}_{m,m'}$ contains the highest weight vector
$z_1^mw_n^{m'}$ which is of the weight $(m,0,\cdots ,0,-m')$.
Therefore, $L_{m,m'}$ contains an irreducible representation
${\hat T}_{m,m'}$ of $U_q({\rm gl}_n)$ with highest weight
$(m,0,\cdots , 0,-m')$. This irreducible representation is absent
in the decomposition (25). Hence, ${\hat T}_{m,m'}$ is a
subrepresentation in $T_{m,m'}$. By the formula for dimensions of
irreducible representations of $U_q({\rm gl}_n)$ and by Corollary
3 we have $\dim {\hat T}_{m,m'}=\dim {\cal H}_{m,m'}$. Therefore,
${\hat T}_{m,m'}$ is equivalent to $T_{m,m'}$. Proposition is
proved.

Thus, we proved that the action of the algebra $U_q({\rm gl}_n)$
on the space ${\cal A}$ realizes the irreducible representations
$T_{m,m'}$ on the subspaces ${\cal H}_{m,m'}$ of homogeneous
$q$-harmonic polynomials, respectively.

We denote by ${\cal A}^{U_q({\rm gl}_n)}$ the space of elements of
${\cal A}$ consisting of invariant elements with respect to the
action of $U_q({\rm gl}_n)$.

{\it Proposition 5:} {\it We have ${\cal A}^{U_q({\rm
gl}_n)}={\Bbb C}[Q]$ and}
\[
{\cal A}^{U_q({\rm gl}_{n-1})}\simeq \bigoplus_{k,l} {\Bbb
C}[Q_{n-1}] z_n^k w_n^l \simeq \bigoplus_{k,l} {\Bbb C}[Q] z_n^k
w_n^l.
\]

{\it Proof:} The formula (23) leads to the decomposition of the
representation $L$ on $\cal A$ into irreducible subrepresentations
of ${U_q({\rm gl}_n)}$ (the representation multiple to the
irreducible representation $T_{m,m'}$ is realized on ${\Bbb
C}[Q]\otimes {\cal H}_{m,m'}$). Since the trivial representation
of ${U_q({\rm gl}_n)}$ is realized only on ${\cal H}_{0,0}$, then
${\cal A}^{U_q({\rm gl}_n)}$ coincides with ${\Bbb C}[Q]\otimes
{\cal H}_{0,0}\equiv$ ${\Bbb C}[Q]\otimes {\Bbb C} \simeq$ ${\Bbb
C}[Q]$.

In order to prove the second equality we note that for
$U_q({\rm gl}_{n-1})$-module ${\cal A}$ we have
$$
{\cal A}={\Bbb C}_q[z_1,w_1,\cdots ,z_n,w_n]=\bigoplus _{k,l}
{\Bbb C}_q[z_1,w_1,\cdots ,z_{n-1},w_{n-1}]
z^k_nw^l_n .
$$
The action of the subalgebra $U_q({\rm gl}_{n-1})$ on monomials
$z^k_nw^l_n$ is trivial. Moreover, ${\Bbb C}[z_1,w_1,\cdots
,z_{n-1},w_{n-1}]^{U_q({\rm gl}_{n-1})} ={\Bbb C}[Q_{n-1}]$, where
$Q_{n-1}=z_1w_1+\cdots +z_{n-1}w_{n-1}$. Since $Q=Q_{n-1}+z_nw_n$,
we have
${\cal A}^{U_q({\rm gl}_{n-1})}\simeq \bigoplus_{k,l} {\Bbb
C}[Q_{n-1}] z_n^k w_n^l \simeq \bigoplus_{k,l} {\Bbb C}[Q] z_n^k
w_n^l$.
Proposition is proved.
\bigskip

\noindent {\sf VI. THE DUAL PAIR $(U_q({\rm sl}_2),U_q({\rm
gl}_n))$}
\medskip

The formulas
$$
ke=q^2ek,\ \ \ \ kf=q^{-2}fk,\ \ \ \
ef-fe=\frac{k-k^{-1}}{q-q^{-1}} \eqno (26)
$$
determine the quantum algebra $U_q({\rm sl}_2)$ generated by the
elements $k,k^{-1},e,f$. Let ${\cal L}({\cal A})$ be the space of
linear operators on the algebra ${\cal A}$. It is directly
verified by means of formula (18) that the operators
$$
\omega (k)=q^n\gamma \bar\gamma ,\ \ \ \ \omega (e)=q^{-n+1} {\hat
Q} ,\ \ \ \ \omega (f)=-\Delta_q \eqno (27)
$$
satisfy relations (26). This means that the algebra homomorphism
$\omega : U_q({\rm sl}_2)\to {\cal L}({\cal A})$ uniquely
determined by formulas (27) is a representation of $U_q({\rm
sl}_2)$.

Since the operators $\omega (k)$, $\omega (e)$, $\omega (f)$
commute with the operators $L(X)$, $X\in U_q({\rm gl}_n)$ we can
introduce the representation $\omega \times L$ of the algebra
$U_q({\rm sl}_2)\times U_q({\rm gl}_n)$ on ${\cal A}$, where $L$
is the above defined natural action of $U_q({\rm gl}_n)$ on ${\cal
A}$. This representation is reducible. Let us decompose it into
irreducible constituents.

By (23), we have ${\cal A}=\bigoplus _{m,m'\ge 0} ({\Bbb
C}[Q]\otimes {\cal H}_{m,m'})$. The subspaces ${\Bbb C}[Q]\otimes
{\cal H}_{m,m'}$ are invariant under $U_{q}({\rm sl}_2)\times
U_q({\rm gl}_n)$, since the space ${\Bbb C}[Q]$ is elementwise
invariant under $U_q({\rm gl}_n)$, and for $f\in {\Bbb C}[Q]$ and
$h_{m,m'}\in {\cal H}_{m,m'}$ we have
$$
{\hat Q}(f(Q)\otimes h_{m,m'})=Qf(Q)\otimes h_{m,m'}, \eqno (28)
$$      $$
\Delta_q (Q^r\otimes h_{m,m'})=q^{n-1}[r][r+m+m'+n-1]
Q^{r-1}\otimes h_{m,m'}, \eqno (29)
$$    $$
\gamma \bar\gamma (Q^{r}\otimes h_{m,m'})=q^{2r+m+m'}(Q^{r}\otimes
h_{m,m'}). \eqno (30)
$$
These formulas show that $U_{q}({\rm sl}_2)$ acts on ${\Bbb C}[Q]$
and $U_q({\rm gl}_n)$ acts on ${\cal H}_{m,m'}$. However, this
action of $U_{q}({\rm sl}_2)$ depends on the component ${\cal
H}_{m,m'}$. Taking the basis
$$
| r\rangle :=q^{-r(n-1)} [r+m+m'+n-1]!^{-1}Q^r, \ \ \ \
r=0,1,2,\cdots ,
$$
in the space ${\Bbb C}[Q]$, we find from (28)--(30) that
$$
\omega (k)|r\rangle =q^{2r+m+m'+n}|r\rangle \ \ \ \ \ \omega (f)|
r\rangle =-[r]| r-1\rangle ,
$$    $$
\omega (e) |r\rangle =[r+m+m'+n]|r+1\rangle .
$$
Comparing this representation with the known irreducible
representations of $U_{q}({\rm sl}_2)$ (see, for example, Ref. 11)
we derive that the irreducible representation of $U_{q}({\rm
sl}_2)$ of the discrete series with lowest weight $m+m'+n$ is
realized on the component ${\Bbb C}[Q]$ of the space ${\Bbb
C}[Q]\otimes {\cal H}_{m,m'}$. We denote this representation of
$U_{q}({\rm sl}_2)$ by $D_{m+m'+n}$.

Thus, we have derived that on the subspace ${\Bbb C}[Q]\otimes
{\cal H}_{m,m'}\subset {\cal A}$ the irreducible representation
$D_{m+m'+n}\times T_{m,m'}$ of the algebra $U_{q}({\rm
sl}_2)\times U_q({\rm gl}_n)$ acts. This means that for the
reducible representation $\omega \otimes L$ we have the following
decomposition into irreducible components:
$$
\omega \times L=\bigoplus ^\infty _{m,m'=0} D_{m+m'+n}\times
T_{m,m'},
$$
that is, each irreducible representation of $U_q({\rm gl}_n)$ in
this decomposition determines uniquely the corresponding
irreducible representation of $U_q({\rm sl}_2)$ and vise versa.
This means that $U_q({\rm sl}_2)$ and $U_q({\rm
gl}_n)$ constitute a {\it dual pair} under the action on ${\cal
A}$.
 \bigskip

\noindent {\sf VII. RESTRICTION OF $q$-HARMONIC POLYNOMIALS ONTO
THE QUANTUM SPHERE}
 \medskip

The associative algebra ${\cal F}(S^{\Bbb C}_{q,n-1})$ generated
by the elements $z_1,\cdots ,z_n,w_1,\cdots ,w_n$ satisfying the
relations (1)--(3) and the relation
$$
z_1w_1+z_1w_1+\cdots +z_nw_n=1
$$
is called {\it the algebra of functions on the quantum sphere}
$S^{\Bbb C}_{q,n-1}$ (see Refs. 6, Chap. 11, and 7). It is clear
that the following canonical algebra isomorphism has place:
$$
{\cal F}(S^{\Bbb C}_{q,n-1}) \simeq {\cal A}/{\cal I} ,
$$
where ${\cal I}$ is the two-sided ideal of ${\cal A}$ generated by
the element $Q-1\equiv \sum _i z_iw_i -1$. We denote by $\tau$ the
canonical algebra homomorphism
$$
\tau : {\cal A}\to {\cal A}/{\cal I}\simeq {\cal F}(S^{\Bbb
C}_{q,n-1}) .
$$
This homomorphism is called the {\it restriction} of polynomials
of ${\cal A}$ onto the quantum sphere $S^{\Bbb C}_{q,n-1}$.

{\it Proposition 6:} {\it We have $\tau {\cal H}\simeq {\cal
F}(S^{\Bbb C}_{q,n-1})$. This means that $\tau :{\cal H}\to {\cal
F}(S^{\Bbb C}_{q,n-1})$ is a one-to-one mapping, that is, the
restriction of a $q$-harmonic polynomial to the sphere $S^{\Bbb
C}_{q,n-1}$ determines this polynomial uniquely.}

{\it Proof:} By Theorem 1, we have
${\cal F}(S^{\Bbb C}_{q,n-1})=\tau {\cal A} =\tau ({\Bbb
C}[Q]\otimes {\cal H})=\tau {\cal H}$.
Since $Q$ is invariant with respect to the action of the algebra
$U_q({\rm gl}_n)$, then the ideal ${\cal I}$ is an invariant
subspace under the action of $U_q({\rm gl}_n)$ on ${\cal A}$.
Therefore, an action of $U_q({\rm gl}_n)$ on ${\cal A}/{\cal I}$
is defined. This action coincides with the action in Ref. 8. The
homomorphism $\tau$ intertwines the action of $U_q({\rm gl}_n)$ on
${\cal A}$ and on ${\cal A}/{\cal I}$. Since $\tau {\cal H}_{m,m'}
\ne \{ 0\}$, then the action of $U_q({\rm gl}_n)$ realizes the
same irreducible representation on ${\cal H}_{m,m'}$ and on $\tau
{\cal H}_{m,m'}$. This means that ${\rm dim}\ {\cal H}_{m,m'}=
{\rm dim}\ \tau {\cal H}_{m,m'}$, that is, the mapping $\tau$ is
one-to-one on ${\cal H}_{m,m'}$. Therefore, it is one-to-one on
${\cal H}$. Proposition is proved.

Proposition 6 allows us to determine a scalar product on
${\cal H}$. For this, we use the invariant functional $h$ on the
quantum sphere defined in Ref. 8, section 4.2. This functional $h$
is determined by introducing a linear gradation in $\tau {\cal
A}$:
$\tau {\cal A} =\sum _{\lambda \in {\Bbb Z}^n} (\tau {\cal
A})^\lambda$,
where $ (\tau {\cal A})^\lambda =\{ p\in \tau {\cal A}\ |\ p({\bf
t}{\bf z}, {\bf t}^{-1}{\bf w})={\bf t}^{\lambda}p ({\bf z},{\bf
w})\}$, ${\bf t}=(t_1,t_2,\cdots ,t_n)$ are $n$ indeterminates,
and
$$
{\bf t}{\bf z} =(t_1z_1,\cdots ,t_nz_n),\ \ \ {\bf t}^{-1}{\bf w}
=(t^{-1}_1w_1,\cdots ,t^{-1}_nw_n),\ \ \ {\bf t}^{\lambda}
=t^{\lambda _1}_1\cdots t^{\lambda _n}_n.
$$
The subalgebra $(\tau {\cal A})^0$ is spanned by the monomials
$z_1^{\mu _1}\cdots z_n^{\mu _n}w_n^{\mu _n}\cdots w_1^{\mu _1}$
(or by the monomials $w_1^{\mu _1}\cdots w_n^{\mu _n}z_n^{\mu
_n}\cdots z_1^{\mu _1}$), $\mu _i=0,1,2,\cdots $. The functional
$h$ is defined as a linear mapping $h: \tau {\cal A}\to {\Bbb C}$
such that $h(p)=0$ if $p\in ((\tau\cal A)^\lambda$, $\lambda\ne
0$, and
$$
h(w_1^{\mu _1}\cdots w_n^{\mu _n}z_n^{\mu _n}\cdots z_1^{\mu _1})
=\frac{(q^2;q^2)_{\mu _1}\cdots (q^2;q^2)_{\mu _n}(q^2;q^2)_{n-1}}
{(q^2;q^2)_{\mu _1+\cdots +\mu _n+n-1}} .
$$
The following assertions are proved in Ref. 8:

(a) The subalgebra $(\tau {\cal A})^0$ is a commutative algebra
generated by the elements $Q_{n-1}$, $Q_{n-2},\cdots , Q_1$.

(b) The algebra $(\tau {\cal A})^0$ is isomorphic to the
polynomial algebra in $n-1$ commuting indeterminates.

(c) For any polynomial $p({\bf z},{\bf w})=f(Q_1,\cdots
,Q_{n-1})\in (\tau {\cal A})^0$ the value $h(p)$ is expressed in
term of Jackson integral:
$$
h(p)=\frac{(q^2;q^2)_{n-1}}{(1-q^2)^{n-1}} \int _0^1 \int
_0^{Q_{n-1}} \cdots \int _0^{Q_2} f(Q_1,\cdots ,Q_{n-1})d_{q^2}Q_1
\cdots d_{q^2}Q_{n-1}
$$
(the definition of Jackson integral see, for example, in Ref. 12,
Chap. 1).

Now we can introduce a scalar product $\langle \cdot ,\cdot
\rangle$ on ${\cal H}$:
$$
\langle p_1,p_2\rangle =h((\tau p_1)(\tau p_2)^*), \eqno (31)
$$
where $a^*$ determines an element conjugate to $a\in {\cal A}$
under action of the $*$-operation.

{\it Proposition 7:} {\it We have ${\cal H}_{m,m'}\bot {\cal
H}_{r,r'}$ if $(m,m')\ne (r,r')$.}

{\it Proof} follows from the fact that $(\tau p_1)(\tau
p_2)^*\not\in (\tau {\cal A})^0$ if $p_1\in {\cal H}_{m,m'}$,
$p_2\in {\cal H}_{r,r'}$, and $(m,m')\ne (r,r')$.
 \bigskip

\noindent {\sf VIII. THE PROJECTION ${\cal A}_{m,m'}\to {\cal
H}_{m,m'}$}
 \medskip

Let us go back to the decomposition (20) and construct the
projector
$$
{\sf H}_{m,m'} :{\cal A}_{m,m'} ={\cal H}_{m,m'}\oplus Q{\cal
A}_{m-1,m'-1}\to {\cal H}_{m,m'}.
$$
We present this projector in the form
$$
{\sf H}_{m,m'} p=\sum _{k=0}^{\min (m,m') } \alpha _k {\hat
Q}^k\Delta_q^k p, \ \ \ \  \alpha _k \in {\Bbb C},\ \ \ p\in {\cal
A}_{m,m'}. \eqno (32)
$$
We have to calculate values of the coefficients $\alpha _k$. In
order to do this, we act by the operator $\Delta_q$ upon both
parts of (32) and use the relation (18). Under this action, the
left hand side vanishes. Equating the right hand side to 0, we
derive a recurrence relation
$$
q^{n-1}[k][m+m'+n-k-1] \alpha _k+\alpha _{k-1}=0
$$
for $\alpha _k$ which gives
$$
\alpha _k =(-1)^kq^{-(n-1)k}\frac{[m+m'+n-k-2]!}{[k]![m+m'+n-2]!},
\eqno (33)
$$
where $[s]!=[s][s-1][s-2]\cdots [1]$ for $s\ne 0$ and $[0]!=1$.

Note that the coefficients $\alpha _k$ are determined by the
recurrence relation uniquely up to a constant. In (33) we have
chosen this constant in such a way that ${\sf H}_{m,m'}p=p$ for
$p\in {\cal H}_{m,m'}$. This means that ${\sf H}_{m,m'}^2={\sf
H}_{m,m'}$.

{\it Proposition 8:} {\it The operator ${\sf H}_{m,m'}$ commutes
with the action of $U_q({\rm gl}_n)$, that is, with the operators
of the representation $L_{m,m'}$ of $U_q({\rm gl}_n)$.}

{\it Proof:} This assertion follows from the fact that the
operators $L_{m,m'}(X)$, $X\in U_q({\rm gl}_n)$, commute with
${\hat Q}$ and $\Delta_q$ (see Proposition~2). Proposition is
proved.

A polynomial $\varphi$ of the space ${\cal H}_{m,m'}$ is called
{\it zonal} if it is invariant with respect to operators
$L_{m,m'}(X)$, $X\in U_q({\rm gl}_{n-1})$. We shall show below
that zonal polynomials can be expressed in terms of the basic
hypergeometric function ${}_2 \varphi _1$ which is defined by the
formula
$$
{}_2 \varphi _1 (a,\ b;\ c;\ q,\ x)=\sum _{k=0}^\infty
{(a;q)_k(b;q)_k\over (c;q)_k(q;q)_k} x^k
$$
(see Refs. 12 and 13 for properties of this function).

{\it Proposition 9:} (a) {\it The subspace of zonal polynomials in
${\cal H}_{m,m'}$ is one-dimensional.} (b) {\it Up to a constant,
a zonal polynomial of ${\cal H}_{m,m'}$ is given by the formula
$$
\varphi '_{m,m'}=z_n^{m-m'}Q^{m'} \sum _{s=0}^{m'}
{(q^{-2m'};q^2)_s(q^{2(m+n-1)};q^2)_s\over (q^{2(n-1)};q^2)_s
(q^2;q^2)_s} {Q_{n-1}^s\over Q^s}q^{2s}    \eqno (34)
$$
if $m\ge m'$ and by the formula
$$
\varphi '_{m,m'}=Q^{m} \sum _{s=0}^{m}
{(q^{-2m};q^2)_s(q^{2(m'+n-1)};q^2)_s\over (q^{2(n-1)};q^2)_s
(q^2;q^2)_s} {Q_{n-1}^s\over Q^s}q^{2s} w_n^{m'-m} ,  \eqno (35)
$$
if $m\le m'$}.

{\it Proof:} (a) As we have seen, the irreducible representation
$T_{m,m'}$ of $U_q({\rm gl}_n)$ with highest weight $(m,0,\cdots
,0,-m')$ is realized on ${\cal H}_{m,m'}$. It is known that this
representation, under restriction to $U_q({\rm gl}_{n-1})$,
contains trivial (one-dimensional) representation of this
subalgebra with multiplicity 1. This proves the first assertion.

(b) We construct a zonal polynomial of ${\cal H}_{m,m'}$ by using
the projection operator ${\sf H}_{m,m'}$. In order to do this, we
have to take a polynomial $p\in {\cal A}_{m,m'}$ invariant with
respect to $U_q({\rm gl}_{n-1})$ and to act upon it by the
operator ${\sf H}_{m,m'}$. Since the projector ${\sf H}_{m,m'}$
commutes with the action of $U_q({\rm gl}_{n-1})$, a polynomial
obtained in this way is a zonal polynomial. Clearly, the
polynomial $p=z_n^mw_n^{m'}$ belongs to ${\cal A}_{m,m'}$ and is
invariant under the action of $U_q({\rm gl}_{n-1})$. In order to
find an expression for ${\sf H}_{mm'}(z_n^mw_n^{m'})$ we first
assume that $m\ge m'$.

Using the second expression for $\Delta_q$ in (17) and relation
$\bar \partial_n \hat z_n= \hat z_n \bar \partial_n$  we have
$$
\varphi _{m,m'} :={\sf H}_{m,m'}(z_n^mw_n^{m'} ) = \sum
_{s=0}^{m'} \alpha _s {\hat Q}^s\Delta_q^s z_n^mw_n^{m'}
$$   $$
=z_n^{m-m'}\sum _{s=0}^{m'} \alpha _s q^{2(n-1)s} {\hat Q}^s {[m]!\over
[m-s]!} {[m']!\over [m'-s]!} z_n^{m'-s}w_n^{m'-s} .
$$
Taking into account the expression for the coefficients $\alpha
_s$ and using the formulas
$$
[s]!={(q^2;q^2)_s(-1)^s\over (q-q^{-1})^s} q^{-s(s+1)/2},
\ \ \ \ \
{[m]!\over [m-s]!} ={(q^{-2m};q^2)_s\over (q-q^{-1})^s} q^{ms-s(s-1)/2}
$$
we obtain
$$
\varphi _{m,m'} =z_n^{m-m'} \sum _{s=0}^{m'}
q^{2s}{(q^{-2m};q^2)_s (q^{-2m'};q^2)_s \over (q^2;q^2)_s
(q^{-2(m+m'+n-2)};q^2)_s} Q^s z_n^{m'-s}w_n^{m'-s} . \eqno (36)
$$
Using the first relation in (16) we obtain from (36) that
$$
\varphi _{m,m'} =Q^{m'}z_n^{m-m'} \sum _{s=0}^{m'}
q^{2s}{(q^{-2m};q^2)_s (q^{-2m'};q^2)_s \over (q^2;q^2)_s
(q^{-2(m+m'+n-2)};q^2)_s} (Q_{n-1}/Q;q^{-2})_{m'-s} .
$$
Since (see relation (II.4) from Appendix II in Ref. 12)
$$
(Q_{n-1}/Q;q^{-2})_{m'-s} =\sum _{\nu =0}^{m'-s} q^{2\nu}
{(q^{-2(m'-s)};q^2)_\nu \over (q^2;q^2)_\nu  } Q_{n-1}^\nu
/Q^\nu ,
$$
we have
$$
\varphi _{m,m'} =Q^{m'}z_n^{m-m'} \sum _{s=0}^{m'}
{q^{2s}(q^{-2m};q^2)_s (q^{-2m'};q^2)_s \over (q^2;q^2)_s
(q^{-2(m+m'+n-2)};q^2)_s } \sum _{\nu =0}^{m'-s} q^{2\nu}
{(q^{-2(m'-s)};q^2)_\nu \over (q^2;q^2)_\nu  } {Q_{n-1}^\nu \over
Q^\nu}
$$    $$
=Q^{m'}z_n^{m-m'} \sum _{\nu =0}^{m'} \frac{Q_{n-1}^\nu}{Q^\nu}
q^{2\nu} \sum _{s=0}^{m'-\nu} {(q^{-2(m'-s)};q^2)_\nu \over
(q^2;q^2)_\nu  } {q^{2s}(q^{-2m};q^2)_s (q^{-2m'};q^2)_s \over
(q^2;q^2)_s  (q^{-2(m+m'+n-2)};q^2)_s} .  \eqno (37)
$$
Applying relation (I.7) and then relation (I.13) from Appendix I
in Ref. 12 we find
$$
(q^{-2(m'-s)};q^2)_\nu =(-1)^\nu q^{-2m'\nu}q^{\nu (\nu -1)}
{(q^{2m'-2\nu +2};q^2)_\nu (q^{-2m'+2\nu};q^2)_s \over
(q^{-2m'};q^2)_s} .
$$
Therefore, for the sum over $s$ in (37) (which will be denoted by
$I_\nu$) we obtain the expression
$$
I_\nu = (-1)^\nu q^{-2m'\nu}q^{\nu (\nu -1)} {(q^{2m'-2\nu
+2};q^2)_\nu \over (q^2;q^2)_\nu } \sum _{s=0}^{m'-\nu}
{(q^{-2m'+2\nu };q^2)_s (q^{-2m};q^2)_s q^{2s} \over (q^2;q^2)_s
(q^{-2(m+m'+n-2)};q^2)_s } .
$$
The sum over $s$ here is the basic hypergeometric function
$$
{}_2 \varphi _1 (q^{-2m},\ q^{-2m'+2\nu};\ q^{-2(m+m'+n-2)};\
q^2;\ q^2 )= {q^{-2mm'+2m\nu}(q^{-2m'-2n +4};q^2)_{m'-\nu}  \over
(q^{-2(m+m'+n-2)};q^2)_{m'-\nu} },
$$
where we used formula (II.6) from Appendix II in Ref. 12.

Therefore, for the function $\varphi _{m,m'}$ we have the
expression
$$
\varphi _{m,m'}=Q^{m'}z_n^{m-m'} \sum _{\nu =0}^{m'}
\frac{Q_{n-1}^\nu}{Q^\nu}  q^{2\nu}(-1)^\nu q^{-2m'\nu}q^{\nu (\nu
-1)} q^{-2mm'+2m\nu}
$$     $$
\times {(q^{2m'-2\nu +2};q^2)_\nu \over (q^2;q^2)_\nu  }
{(q^{-2m'-2n +4};q^2)_{m'-\nu}  \over
(q^{-2(m+m'+n-2)};q^2)_{m'-\nu} }.
$$
By formula (I.8) of Appendix I in Ref. 12 we have
$$
(q^{2m'-2\nu +2};q^2)_\nu =(q^{-2m'};q^2)_\nu (-1)^\nu q^{2m'\nu}
q^{-\nu (\nu -1)}
$$
and by formula (I.11) from Appendix I in Ref. 12 we obtain
$$
{(q^{-2m'-2n +4};q^2)_{m'-\nu}  \over
(q^{-2(m+m'+n-2)};q^2)_{m'-\nu} }=q^{-2m\nu} {(q^{-2m'-2n
+4};q^2)_{m'} \over (q^{-2(m+m'+n-2)};q^2)_{m'} }
{(q^{2(m+n-1)};q^2)_{\nu}  \over (q^{2n -2};q^2)_{\nu} } .
$$
For this reason, we have
$$
{\sf H}_{mm'}(z_n^mw_n^{m'})= \varphi _{m,m'}= q^{-2mm'}
{(q^{-2m'-2n +4};q^2)_{m'}  \over (q^{-2(m+m'+n-2)};q^2)_{m'} }
Q^{m'}z_n^{m-m'} \times
$$     $$
\times \sum _{\nu =0}^{m'} {(q^{-2m'};q^2)_\nu (q^{2(m+n
-1)};q^2)_{\nu} \over (q^2;q^2)_\nu (q^{2(n -1)};q^2)_{\nu} }
\frac{Q_{n-1}^\nu}{Q^\nu}  q^{2\nu}
$$  $$
=\frac{(q^{2(n-1)};q^2)_{m'}} {(q^{2(m+n-1)};q^2)_{m'}}
Q^{m'}z_n^{m-m'} {}_2\varphi_1(q^{-2m'},q^{2(m+n-1)};q^{2(n-1)};
q^2;q^2Q_{n-1}/Q).
$$

This proves the second assertion of the proposition for the case
$m\ge m'$. The case $m<m'$ is proved in the same way. Proposition
is proved.

The formula
$$
P^{(\alpha ,\beta )}_k(x;\ q)= {}_2 \varphi _1 (q^{-k},\ q^{\alpha
+\beta +k+1};\ q^{\alpha +1};\ q;\ qx)
$$
defines the so called little $q$-Jacobi polynomials. The zonal
polynomials from Proposition 9 can be written in term of these
polynomials as
$$
\varphi '_{m,m'}=Q^{m'}z_n^{m-m'}P^{( n-2,m-m')}_{m'}(Q_{n-1}/Q;\ q^2)
$$
if $m\ge m'$ and as
$$
\varphi '_{m,m'}=Q^{m}P^{( n-2,m'-m)}_{m}(Q_{n-1}/Q;\ q^2)w_n^{m'-m}
$$
if $m\le m'$. Restricting these polynomials onto the quantum
sphere $S^{\Bbb C}_{q,n-1}$ we obtain
$$
\tau \varphi '_{m,m'}= z_n^{m-m'} P^{( n-2,m-m')}_{m'}(Q_{n-1};\
q^2)
$$
if $m\ge m'$ and as
$$
\tau \varphi '_{m,m'}=P^{( n-2,m'-m)}_{m}(Q_{n-1};\ q^2)w_n^{m'-m}
$$
if $m\le m'$. These polynomials are called {\it zonal spherical
functions} on the quantum sphere $S^{\Bbb C}_{q,n-1}$ and were
calculated in Ref. 8 (see also Refs. 3 and 4).
 \bigskip

\noindent {\sf IX. $q$-ANALOGUE OF ASSOCIATED SPHERICAL HARMONICS
WITH RESPECT TO $U_q({\rm gl}_{n-1})$}
 \medskip

It is known (see Ref. 1, Chap. 11) that in the space of classical
homogeneous harmonic polynomials on the unitary (complex
Euclidean) space $E_n^{\Bbb C}$ there exist different orthonormal
bases. They correspond to different separations of variables. Each
separation of variables corresponds to a certain chain of
subgroups of the unitary group $U(n)$. A similar picture has place
for the spaces ${\cal H}_{m,m'}$ of homogeneous $q$-harmonic
polynomials. We consider in this section a $q$-analogue of
separation of variables corresponding to spherical coordinates on
the sphere $S_{n-1}^{\Bbb C}$ (see Ref. 1, Chap. 11).

In the classical case, the tree method distinguishes different
separations of variables. Different separations of variables are
in a one-to-one correspondence with different chains of
subgroups of $U(n)$. The same tree method can be used for
$q$-harmonic polynomials, but instead of chains of subgroups of
$U(n)$ we have to take the corresponding chains of subalgebras of
the algebra $U_q({\rm gl}_{n})$. A certain orthogonal basis
corresponds to such a chain of subalgebras.

The aim of this section is to construct an orthogonal basis of the
space ${\cal H}_{m,m'}$ of homogeneous $q$-harmonic polynomials
which corresponds to the chain
$$
U_q({\rm gl}_n)\supset U_q({\rm gl}_{n-1})\supset \cdots \supset
U_q({\rm gl}_3) \supset U_q({\rm gl}_2)\supset U_q({\rm gl}_1).
\eqno (38)
$$
This basis is a $q$-analogue of the set of associated spherical
harmonics on the complex vector space which are products of
certain Jacobi polynomials (see, Ref. 1, Chap. 11).
The basis elements give solutions of the equation $\Delta _qp=0$
in "separated coordinates". So, we obtain a $q$-analogue of the
classical separation of variables.

{\it Lemma 1:} {\it Let $f_{k'}({\bf z}')$ and $g_{l'}({\bf w}')$
be homogeneous polynomials of degrees $k'$ in ${\bf z}'\equiv
(z_1,z_1,\cdots,$ $z_{n-1})$ and of degrees $l'$ in ${\bf
w}'\equiv (w_1,w_1,\cdots ,w_{n-1})$, respectively. Then for any
nonnegative integers $k$ and $l$ we have
$$
\Delta_q (z_n^kw_n^l f_{k'}({\bf z}')g_{l'}({\bf w}') )= q^{l-k}
z_n^kw_n^l \Delta_{n-1} (f_{k'}({\bf z}')g_{l'}({\bf w}') )
$$  $$
+q^{2(n-1)}q^{l'+k'}[k][l]z_n^{k-1}w_n^{l-1} f_{k'}({\bf z}')
g_{l'}({\bf w}'),
$$
where $\Delta _{n-1}=\sum _{i=1}^{n-1}q^{2i-2}{\bar\partial}_i
\partial_i$ is
the $q$-Laplace operator for the elements ${\bf z}'\equiv
(z_1,$ $\cdots,$ $z_{n-1})$ and ${\bf w}'\equiv (w_1,\cdots
,w_{n-1})$.}

{\it Proof:} Using the relations for the operators from section III
we derive
\[
\bar \partial_n \partial_n z_n^kw_n^l f_{k'}({\bf z}')g_{l'}({\bf
w}')= q^{(l-k)k'} \bar \partial_n  \partial_n  f_{k'}({\bf
z}')z_n^k w_n^l g_{l'}({\bf w}')
\]
\[
=q^{(l-k)k'}[k] \bar \partial_n f_{k'}({\bf z}')z_n^{k-1} w_n^l
g_{l'}({\bf w}')= q^{(l-1)k'}[k] \bar \partial_n  z_n^{k-1}
f_{k'}({\bf z}') w_n^l g_{l'}({\bf w}')
\]
\[
=q^{(l-1)k'}[k]  z_n^{k-1} \bar \partial_n  f_{k'}({\bf z}') w_n^l
g_{l'}({\bf w}')= q^{-2k'+l(k'+l')}[k]  z_n^{k-1}
  f_{k'}({\bf z}') \bar \partial_n g_{l'}({\bf w}') w_n^l
\]\[
=q^{-2k'+l k'+l'}[k][l]  z_n^{k-1}
  f_{k'}({\bf z}') w_n^{l-1} g_{l'}({\bf w}') =
q^{-k'+l'}[k][l]  z_n^{k-1} w_n^{l-1}
  f_{k'}({\bf z}') g_{l'}({\bf w}') .
\]
Since $\partial_i \hat w_n = q \hat w_n
\partial_i + (1-q^2) \hat w_i\partial_n$, $i<n$, and
$\partial_n(w_n^l f_{k'}({\bf z}')g_{l'}({\bf w}'))=0$, we have
\[
\bar \partial_i  \partial_i z_n^kw_n^l f_{k'}({\bf z}')g_{l'}({\bf
w}')= \bar \partial_i   z_n^k \partial_iw_n^l f_{k'}({\bf
z}')g_{l'}({\bf w}')= q^l \bar \partial_i  z_n^k w_n^l \partial_i
f_{k'}({\bf z}')g_{l'}({\bf w}').
\]
Using reccurently the relation
$\bar \partial_i \hat z_n= q^{-1} \hat z_n \bar \partial_i
+(1-q^{-2}) q^{2(n-i)} \hat z_i\bar \partial_n$,
we obtain
\[
\bar \partial_i  \partial_i z_n^kw_n^l f_{k'}({\bf z}')g_{l'}({\bf w}')=
q^{l-k}
z_n^kw_n^l \bar\partial_i  \partial_i f_{k'}({\bf z}')g_{l'}({\bf w}')
\]
\[
+q^{2(n-i)}(1-q^{-2}) q^{l'+1}[k][l] z_n^{k-1}w_n^{l-1}\hat z_i\partial_i
 f_{k'}({\bf z}')g_{l'}({\bf w}').
\]
Thus, one has
\[
\Delta_{n-1} (z_n^kw_n^l f_{k'}({\bf z}')g_{l'}({\bf w}'))=
q^{l-k} z_n^kw_n^l \Delta _{n-1} f_{k'}({\bf
z}')g_{l'}({\bf w}')
\]
\[
+q^{2n-3} (q^2-1) q^{l'} [k][l][k'] z_n^{k-1} w_n^{l-1}
f_{k'}({\bf z}')g_{l'}({\bf w}'),
\]
where the relation
\[
\sum_{i=1}^{n-1} \hat z_i \partial_i f_{k'}({\bf z}')g_{l'}({\bf
w}')= \sum_{i=1}^{n} \hat z_i \partial_i f_{k'}({\bf z}')g_{l'}({\bf
w}')=
\{\gamma \} f_{k'}({\bf z}')g_{l'}({\bf w}')= [k']
f_{k'}({\bf z}')g_{l'}({\bf w}')
\]
has been used. From the above results and from the equality
$\Delta_q= q^{2(n-1)}\bar\partial_n  \partial_n + \Delta_{n-1}$,
the lemma follows.

{\it Proposition 10:} {\it Let $s$ and $s'$ be integers such that
$0\le s\le m$ and $0\le s'\le m'$. Let $h_{s,s'}({\bf z}',{\bf
w}')$ be a homogeneous harmonic polynomial of degree $s$ in ${\bf
z}'=(z_1,z_2,\cdots ,z_{n-1})$ and of degree $s'$ in ${\bf
w}'=(w_1,w_2,\cdots ,w_{n-1})$. Then for $z^{m-s}_n
w^{m'-s'}_nh_{s,s'}({\bf z}',{\bf w}')\in {\cal A}_{m,m'}$ we have
$$
{\sf H}_{m,m'}(z^{m-s}_n w^{m'-s'}_nh_{s,s'}({\bf z}',{\bf w}'))=
z_n^{m-s-m'+s'}Q^{m'-s'}d^{mm'}_{ss'}h_{s,s'}({\bf z}',{\bf w}')
, \eqno (39)
$$
where $m-s\ge m'-s'$,
$$
d^{mm'}_{ss'}=q^{-2(m-s)(m'-s')}\frac{(q^{-2m'-2s-2n+4};q^2)_{m'-s'}}
{(q^{-2m-2m'-2n+4};q^2)_{m'-s'}} \times
$$      $$
\times {}_2\varphi _1(q^{-2(m'-s')},q^{2(m+s'+n-1)};
q^{2(s+s'+n-1)};\ q^2;\ q^2Q_{n-1}/Q),
$$
and
$$
{\sf H}_{m,m'}(z^{m-s}_n w^{m'-s'}_nh_{s,s'}({\bf z}',{\bf w}'))=
Q^{m-s}d^{mm'}_{ss'}w_n^{m'-s'-m+s}h_{s,s'}({\bf z}',{\bf w}') ,
\eqno (40)
$$
where $m-s\le m'-s'$,}
$$
d^{mm'}_{ss'}=q^{-2(m-s)(m'-s')}\frac{(q^{-2m-2s'-2n+4};q^2)_{m-s}}
{(q^{-2m-2m'-2n+4};q^2)_{m-s}} \times
$$      $$
\times {}_2\varphi _1(q^{-2(m-s)},q^{2(m'+s+n-1)};
q^{2(s+s'+n-1)};\ q^2;\ q^2Q_{n-1}/Q).
$$

{\it Proof:} The proof of this proposition is similar to that of
Proposition 6 and we shall omit details. Taking into account
formula (32) for the projector ${\sf H}_{m,m'}$ and Lemma 1, we
obtain
$$
{\sf H}_{m,m'}(z^{m-s}_n w^{m'-s'}_nh_{s,s'}({\bf z}',{\bf w}'))=
\sum_{k=0}^{\min(m,m')} \alpha_k Q^k \Delta_q^k z^{m-s}_n
w^{m'-s'}_nh_{s,s'}({\bf z}',{\bf w}')
$$ $$
=\sum_{k=0}^b \alpha_k Q^k q^{2(n-1)k} q^{(s+s')k}
\frac{[m-s]![m'-s']!}{[m-s-k]![m'-s'-k]!} z^{m-s-k}_n
w^{m'-s'-k}_nh_{s,s'}({\bf z}',{\bf w}'),
$$
where $b=\min(m-s,m'-s')$. Let $m-s\ge m'-s'$, then
$$
{\sf H}_{m,m'}(z^{m-s}_n w^{m'-s'}_nh_{s,s'}({\bf z}',{\bf w}'))=
z_n^{m-s-m'+s'}Q^{m'-s'}d^{mm'}_{ss'}h_{s,s'}({\bf z}',{\bf w}'),
$$
where
$$
d^{mm'}_{ss'}=\sum _{k=0}^{m'-s'} q^{2k}\frac{(q^{-2(m-s)};q^2)_k
(q^{-2(m'-s')};q^2)_k} {(q^2;q^2)_k (q^{-2(m+m'+n-2)};q^2)_k} \sum
_{\nu =0}^{d} q^{2\nu} \frac{(q^{-2(m'-s'-k)};q^2)_\nu }
{(q^2;q^2)_\nu } \frac{Q^\nu_{n-1}}{Q^\nu} .
$$
Here $d=m'-s'-k$. Changing the order of summations in the last
expression we have
$$
d^{mm'}_{ss'}=\sum _{\nu =0}^{\sigma '}\frac{(Q_{n-1}/Q)^\nu
q^{2\nu}} {(q^2;q^2)_\nu } \sum _{k=0}^{\sigma '-\nu} q^{2k}
$$      $$
\times \frac{(q^{-2(m-s)};q^2)_k  (q^{-2(m'-s'-k)};q^2)_\nu
(q^{-2(m'-s')};q^2)_k } {(q^2;q^2)_k (q^{-2(m+m'+n-2)};q^2)_k} ,
\eqno (41)
$$
where $\sigma '=m'-s'$. Since
$$
(q^{-2(m'-s'-k)};q^2)_\nu =q^{\nu (\nu -1)}(-q^{-2(m'-s'-k)})^\nu
(q^{2m'-2s'-2k}q^{-2\nu +2};q^2)_\nu
$$          $$
=(-1)^\nu q^{\nu (\nu -1)}q^{-2(m'-s')\nu} \frac{(q^{2m'-2s'-2\nu
+2};q^2)_\nu (q^{-2m'+2s'+2\nu };q^2)_k} {(q^{-2m'+2s'};q^2)_k} ,
$$
for the sum over $k$ in (41) we have
$$
\frac{(-1)^\nu q^{\nu (\nu -1)}
(q^{-2\nu}q^{2m'-2s'+2};q^2)_{\nu}}{q^{2(m'-s')\nu}}
\sum _{k=0}^{\sigma '-\nu} \frac{(q^{-2m'+2s'+2\nu };q^2)_k
(q^{-2(m-s) };q^2)_k}{(q^2;q^2)_k (q^{-2(m+m'+n-2)};q^2)_k} q^{2k}
$$
$$
=a_\nu(-1)^\nu q^{\nu (\nu -1)}q^{-2(m'-s')\nu}
{}_2 \varphi _1(q^{-2(m-s)}, q^{-2(m'-s')+2\nu};q^{-2(m+m'+n-2)};\ q^2;\
q^2)
$$        $$
=a_\nu(-1)^\nu q^{\nu (\nu -1)}q^{-2(m'-s')\nu}
\frac{ (q^{-2m'-2s-2n+4 };q^2)_{m'-s'-\nu}}
{ (q^{-2m-2m'-2n+4 };q^2)_{m'-s'-\nu}} q^{-2(m-s)(m'-s'-\nu)}
$$       $$
=a_\nu\frac{(-1)^\nu q^{\nu (\nu -1)}}{q^{2(m'-s')\nu}} \frac{
(q^{-2m'-2s-2n+4 };q^2)_{m'-s'}} { (q^{-2m-2m'-2n+4
};q^2)_{m'-s'}} \frac{ (q^{2m+2s'+2n-2 };q^2)_{\nu}} {
(q^{2s+2s'+2n-2 };q^2)_{\nu}} q^{-2(m-s)(m'-s')},
$$
where $a_\nu =(q^{-2\nu}q^{2m'-2s'+2};q^2)_{\nu}$. Since
$$
a_\nu \equiv (q^{-2\nu}q^{2m'-2s'+2};q^2)_{\nu} =
(-1)^\nu q^{-\nu (\nu -1)}q^{2(m'-s')\nu} (q^{-2m'+2s'};q^2)_{\nu} ,
$$
for $d^{mm'}_{ss'}$ we have the expression
$$
d^{mm'}_{ss'}=c^{mm'}_{ss'}
{}_2 \varphi _1 (q^{-2(m'-s')}, q^{2(m+s'+n-1)};
q^{2(s+s'+n-1)};\ q^2;\ q^2Q_{n-1}/Q),
$$
where
$$
c^{mm'}_{ss'}=
q^{-2(m-s)(m'-s')}
\frac{ (q^{-2m'-2s-2n+4 };q^2)_{m'-s'}}
{ (q^{-2m-2m'-2n+4 };q^2)_{m'-s'}}=
\frac{ (q^{2(s+n-1)} ;q^2)_{m'-s'}}
{ (q^{2(m+n-1) };q^2)_{m'-s'}}.
$$
In the case when $m-s\le m'-s'$, the proof is similar and we omit
it. Proposition is proved.

{\it Remark:} If $n=2$, then polynomials $h_{s,s'}(z_1,w_1)$ in
Proposition 10 are multiple to elements from (21), that is, we
have $s=0$ or $s'=0$ or $s=s'=0$.
\medskip

The expressions for $d^{mm'}_{ss'}$ from Proposition 10 can be
represented in terms of little $q$-Jacobi polynomials $P^{(\alpha
,\beta)}_k(x;q)$ as
$$
d^{mm'}_{ss'}=c^{mm'}_{ss'} P^{(s+s'+n-2,m-s-m'+s')}_{m'-s'}(
Q_{n-1}/Q ) ,
$$
if $m-s\ge m'-s'$ and as
$$
d^{mm'}_{ss'}=c^{m'm}_{s's} P^{(s+s'+n-2,m'-s'-m+s)}_{m-s}(
Q_{n-1}/Q ) ,
$$
if $m-s\le m'-s'$.

We denote the expression $z_n^{m-s+m'+s'}Q^{m'-s'}d^{mm'}_{ss'}$
from (39) and the expression $Q^{m-s}d^{mm'}_{ss'}w_n^{m'-s'+m+s}$
from (40) by $t^{n;m,m'}_{s,s'}$. Then
$$
{\sf H}_{m,m'}(z^{m-s}_n w^{m'-s'}_nh_{s,s'}({\bf z}',{\bf w}'))=
t^{n;m,m'}_{s,s'}h_{s,s'}({\bf z}',{\bf w}') . \eqno (42)
$$
Moreover, the space ${\cal H}_{m,m'}$ can be represented as the
direct sum
$$
{\cal H}_{m,m'} =\bigoplus _{s=0}^m \bigoplus _{s'=0}^{m'}
t^{n;m,m'}_{s,s'} {\cal H}^{(n-1)}_{s,s'} , \eqno (43)
$$
where ${\cal H}^{(n-1)}_{s,s'}$ are the corresponding spaces of
homogeneous $q$-harmonic polynomials in $z_i$, $w_i$,
$i=1,2,\cdots ,n-1$. To prove this, we note that the subspaces
$t^{n;m,m'}_{s,s'} {\cal H}^{(n-1)}_{s,s'}$ pairwise do not
intersect and $\bigoplus _{s=0}^m \bigoplus _{s'=0}^{m'}
t^{n;m,m'}_{s,s'} {\cal H}^{(n-1)}_{s,s'}\subset {\cal H}_{mm'}$.
Now the equality (43) follows from the fact that dimensions of the
spaces on the right and on the left coincide.

To have a correspondence with the classical case, below we denote
$t^{2;m,m'}_{s,s'}$ (in this case $s=0$ or $s'=0$) by
$t^{2;m,m'}_{s}$ if $s'=0$ and by $t^{2;m,m'}_{-s'}$ if $s=0$.

Taking into account the orthogonality relation (7.3.3) in Ref. 12
for little $q$-Jacobi polynomials we obtain for the scalar product of
$t^{n;m,m'}_{s,s'}h^{(n-1)}_{s,s'}$ and
$t^{n;m,m'}_{r,r'}h^{(n-1)}_{r,r'}$, $h^{(n-1)}_{p,p'}\in
{\cal H}^{(n-1)}_{p,p'}$, the expression
$$
\langle t^{n;m,m'}_{s,s'}h^{(n-1)}_{s,s'},t^{n;m,m'}_{r,r'}h^{(n-1)}_{r,r'}
\rangle =\delta_{sr}\delta_{s'r'}(c^{mm'}_{ss'})^{-2}b^{mm'}_{ss'}
\langle h^{(n-1)}_{s,s'},h^{(n-1)}_{r,r'}\rangle_{(n-1)},
$$
where $\langle \cdot ,\cdot\rangle_{(n-1)}$ is the scalar product in the
space ${\cal H}^{(n-1)}_{ss'}$ and
$$
b^{mm'}_{ss'}=\frac{(1-q^{2(n+s+s'-1)})q^{2(m'-s')(n+s+s'-1)}
(q^2;q^2)_{m-s}(q^2;q^2)_{m'-s'}}
{(1-q^{2(2m+n-1)})(q^{2(n+s+s'-1)};q^2)_{m-s}
(q^{2(n+s+s'-1)};q^2)_{m'-s'}} .
$$
Note that a calculation of this scalar product reduces to
$q$-integration (see Refs. 3 and 4 on calculation of $q$-integrals of
this type).

Now we apply the decomposition (43) to the subspaces ${\cal
H}^{(n-1)}_{s,s'}$ and obtain
$$
{\cal H}_{m,m'}=\bigoplus _{s=0}^m \bigoplus _{s'=0}^{m'}
\bigoplus _{r=0}^s \bigoplus _{r'=0}^{s'} t^{n;m,m'}_{s,s'}
t^{n-1;s,s'}_{r,r'} {\cal H}^{(n-2)}_{r,r'} ,
$$
where ${\cal H}^{(n-2)}_{r,r'}$ are the subspaces of homogeneous
$q$-harmonic polynomials in $z_i, w_i$, $i=1,2,\cdots ,n-2$.
Continuing such decompositions we obtain the decomposition
$$
{\cal H}_{m,m'}=\bigoplus _{{\bf m},{\bf m}',m_1} {\Bbb C}\Xi
_{{\bf m},{\bf m}',m_1} ({\bf z},{\bf w}),
$$
where the polynomials $\Xi _{{\bf m},{\bf m}',m_1}$ are given by
the formula
$$
\Xi _{{\bf m},{\bf m}',m_1} ({\bf z},{\bf w})=
t^{n;m,m'}_{m_{n-1},m'_{n-1}}
t^{n-1;m_{n-1},m'_{n-1}}_{m_{n-2},m'_{n-2}}\cdots t^{3;m_3,
m'_3}_{m_2,m'_2}t^{2;m_2, m'_2}_{m_1} t^{1;m_1} ,\eqno (44)
$$
and the summation is over all sets of $2n-3$ integers ${\bf
m}=(m_{n-1},\cdots ,m_2)$, ${\bf m}'=(m'_{n-1},\cdots ,m'_2)$,
$m_1$ such that $m_i\ge 0$, $m'_i\ge 0$, $i=2,3,\cdots ,n-1$,
$m_2\ge m_1\ge -m'_2$,
$$
m\ge m_{n-1}\ge m_{n-2}\ge \cdots \ge m_2,\ \ \ m'\ge m'_{n-1}\ge
m'_{n-2}\ge \cdots \ge m'_2.
$$
Here $t^{p;m_{p},m'_{p}}_{m_{p-1},m'_{p-1}}$ and
$t^{2;m_2,m'_2}_{m_1}$ are determined by formulas given above and
$$
t^{1;m_1}=z_1^{m_1}\ \ {\rm for}\ \ m_1>0,\ \ \ \ \ t^{1;0}=1,\ \
\ \ \ t^{1;m_1}=w_1^{-m_1}\ \ {\rm for}\ \ m_1<0.
$$

It is easy to show that the basis (44) is orthogonal with respect
to the scalar product introduced above.

At $q=1$, polynomials (44) turn into the basis elements of the
spaces of homogeneous harmonic polynomials on ${\Bbb C}^n$ in
separated coordinates determined by formulas (2) of section 11.1.4
in Ref. 1. These classical homogeneous harmonic polynomials,
restricted to the sphere $S^{\Bbb C}_{n-1}$, coincide with
associated spherical functions from section 11.3 in Ref. 1. They
are matrix elements of zero column of the corresponding
irreducible representations of the group $U(n)$.

The basis elements (44) give solutions of the equation $\Delta
p=0$ in ${\cal H}_{m,m'}$. A representation of solutions in the
form (44) can be considered as a $q$-analogue of the corresponding
classical separation of variables.

In order to have an orthonormal basis in ${\cal H}_{m,m'}$
we replace each $t^{n-i;m_{n-i},m'_{n-i}}_{m_{n-i-1},m'_{n-i-1}}$
in the expression (44) for
$\Xi _{{\bf m},{\bf m}',m_1} ({\bf z},{\bf w})$ by
$$
{\hat t}^{n-i;m_{n-i},m'_{n-i}}_{m_{n-i-1},m'_{n-i-1}}=
c^{m_{n-i},m'_{n-i}}_{m_{n-i-1},m'_{n-i-1}}
(b^{m_{n-i},m'_{n-i}}_{m_{n-i-1},m'_{n-i-1}})^{-1/2}
t^{n-i;m_{n-i},m'_{n-i}}_{m_{n-i-1},m'_{n-i-1}}.
$$
We denote the expression (44) with such the replacement by
${\hat\Xi}_{{\bf m},{\bf m}',m_1} ({\bf z},{\bf w})$.
These polynomials constitute an orthonormal basis of ${\cal H}_{m,m'}$.

It was shown above that the irreducible representation $T_{m,m'}$
with highest weight $(m,0,\cdots ,0,-m')$ acts on the space ${\cal
H}_{m,m'}$. The following assertion is true.

{\it Proposition 11:} {\it The operators $T_{m,m'}(e_j)$,
$T_{m,m'}(f_j)$ and $T_{m,m'}(k_j)$, corresponding to the generating
elements $e_j, f_j, k_j$ of the algebra $U_q({\rm gl}_n)$,
act upon the basis elements
${\hat\Xi} _{{\bf m},{\bf m}',m_1}\equiv |{\bf m},{\bf m}',m_1\rangle$
as
$$
T_{m,m'}(e_{j-1})|{\bf m},{\bf m}',m_1\rangle =A({\bf m},{\bf m}')
|{\bf m}^{+1}_{j-1},{\bf m}',m_1\rangle +B({\bf m},{\bf m}') |{\bf
m},{{\bf m}'}^{-1}_{j-1},m_1\rangle ,
$$   $$
T_{m,m'}(f_{j-1})|{\bf m},{\bf m}',m_1\rangle =A({\bf
m}^{-1}_{j-1},{\bf m}') |{\bf m}^{-1}_{j-1},{\bf m}',m_1\rangle
+B({\bf m},{{\bf m}'}^{+1}_{j-1}) |{\bf m},{{\bf
m}'}^{+1}_{j-1},m_1\rangle ,
$$  $$
T_{m,m'}(k_{j-1})|{\bf m},{\bf m}',m_1\rangle
=q^{m'_j-m_j+m_{j-1}-m'_{j-1}}|{\bf m},{\bf m}',m_1\rangle ,
$$
where
$$
\leqno A({\bf m},{\bf m}')
$$   $$
=\left( \frac{[m_j-m_{j-1}] [m'_j+m_{j-1}+j-1]
[m_{j-1}-m_{j-2}+1][m_{j-1}+m'_{j-2}+j-2]}{[m_{j-1}+m'_{j-1}+j-2]
[m_{j-1}+m'_{j-1}+j-1]}\right) ^{1/2} ,
$$   $$
\leqno B({\bf m},{\bf m}')
$$   $$
=\left( \frac{[m'_j-m'_{j-1}+1] [m_j+m'_{j-1}+j-2]
[m'_{j-1}-m'_{j-2}][m'_{j-1}+m_{j-2}+j-3]}{[m_{j-1}+m'_{j-1}+j-2]
[m_{j-1}+m'_{j-1}+j-3]}\right) ^{1/2} ,
$$
$m_n\equiv m$, $m'_n\equiv m'$, ${\bf m}_j^{\pm 1}$ denotes the
set of the numbers ${\bf m}_{j-1}$ with $m_{j-1}$ replaced by
$m_{j-1}\pm 1$, respectively.}

A proof of this proposition is awkward. Since it is similar to that
of Theorem 1 in Ref. 5, we omit it.
 \bigskip

\noindent {\sf X. $q$-ANALOGUE OF ASSOCIATED SPHERICAL HARMONICS
WITH RESPECT TO\\ $U_q({\rm gl}_{p}) \times U_q({\rm gl}_{n-p})$}
\medskip

In section IX we found an orthogonal basis of the space ${\cal
H}_{mm'}$ of homogeneous $q$-harmonic polynomials corresponding to
the chain of subalgebras (38). In this section we shall find
orthogonal bases of the same space corresponding to the reductions
$$
U_q({\rm gl}_n)\supset U_q({\rm gl}_{p})\times U_q({\rm
gl}_{n-p})\supset \cdots . \eqno (45)
$$
In the classical case (see Ref. 1, Chap. 11), further
reductions can be made taking any chain of subgroups of the groups
$U(p)$ and $U(n-p)$.
In particular, the usual tree method (see Ref. 1,
section 10.2) can be used to describe different chains of
these groups corresponding to different orthogonal bases of
${\cal H}_{mm'}$.
In our case, there are some difficulties with construction
of orthogonal bases corresponding to any chain of subalgebras in (45).
For this reason, we construct orthogonal bases corresponding to the case,
when we take chains of the type (38) for the subalgebras
$U_q({\rm gl}_{p})$ and $U_q({\rm gl}_{n-p})$ in (45).

We represent the set $({\bf z},{\bf w})=(z_1,\cdots
,z_n;w_1,\cdots ,w_n)$ as $({\bf y},{\bf t})$,
where ${\bf y}=(z_1,z_2$, $\cdots$, $z_p$, $w_1,w_2,\cdots ,w_p)$ and
${\bf t}=(z_{p+1},\cdots ,z_n,w_{p+1},\cdots ,w_n)$. Then the
$q$-Laplace operator $\Delta_q$ can be written as
$$
\Delta_q =\Delta _{({\bf y})} + \Delta _{({\bf t})} , \eqno (46)
$$
where
$$
\Delta _{({\bf y})}=\partial_1 \bar\partial_1+\cdots
+\partial_p\bar\partial_p,\ \ \ \Delta _{({\bf t})}=
\partial_{p+1} \bar\partial_{p+1}+\cdots
+\partial_n\bar\partial_n
=\sum _{i=1}^{n-p}q^{2(i-1)} \bar\partial_{p+i}\partial_{p+i}. \eqno(47)
$$
The operator $\Delta _q$ can be also represented as
$$
\Delta_q =\hat\Delta _{({\bf y})} +q^{2p} \Delta _{({\bf t})} ,
$$
where
$$
\hat\Delta _{({\bf y})}=\bar\partial_1\partial_1+q^2
\bar\partial_2\partial_2+\cdots +q^{2(p-1)}
\bar\partial_{p}\partial_{p}. \eqno(48)
$$
We have
$$
\Delta _{({\bf y})} -\hat\Delta _{({\bf y})}=(1-q^{2p})
\Delta _{({\bf t})}.\eqno(49)
$$

In order to find bases of ${\cal H}_{m,m'}$ corresponding to the
reduction (45) we take nonnegative numbers $r,r',s,s'$ such that
$$
u:=m-r-s=m'-r'-s'\ge 0.
$$
We wish to find a harmonic projection of the polynomials
$$
Q_{\bf y}^u h_{s,s'}({\bf t}) h_{r,r'}({\bf y}) \in {\cal
A}_{mm'}, \ \ \ \  h_{s,s'}({\bf t})\in \tilde{\cal
H}^{({\bf t})}_{ss'},\ \  h_{r,r'}({\bf y}) \in {\cal H}^{({\bf
y})}_{r,r'}, \eqno (50)
$$
where $Q_{\bf y}:=z_1w_1+ \cdots +z_pw_p$, ${\cal H}^{({\bf
y})}_{r,r'}$ is the space of homogeneous $q$-harmonic polynomials
in ${\bf y}=(z_1,z_2,\cdots ,z_p,w_1,w_2,\cdots ,w_p)$, and
$\tilde{\cal H}^{({\bf t})}_{s,s'}$ is the space obtained in the
following way. We take the space ${\cal H}^{(n-p)}_{s,s'}$ of
homogeneous $q$-harmonic polynomials in $(z_1,\cdots ,z_{n-p},
w_1,\cdots ,w_{n-p})$ and, using the relations between
$z_i$ and $w_j$, represent each its polynomial
in such a form that in each of its summands (monomials)
the elements $z_1,\cdots ,z_{n-p}$ stand before the elements
$w_1,\cdots ,w_{n-p}$. Then we replace $z_1,\cdots ,z_{n-p},
w_1$, $\cdots$, $w_{n-p}$ by $z_{p+1},\cdots ,z_{n},
w_{p+1},\cdots ,w_{n}$, respectively, in each of these polynomials.
The space of these polynomials in $z_{p+1},\cdots ,z_{n},
w_{p+1},\cdots ,w_{n}$ is denoted by $\tilde{\cal
H}^{({\bf t})}_{ss'}$.

{\it Lemma 2:} {\it Polynomials $P$ of $\tilde{\cal H}^{({\bf t})}_{s,s'}$
satisfy the conditions $\partial_i P=0$,
$\bar\partial_i P=0$, $i=1,2,\cdots ,p$.}

{\it Proof:} Fulfillment of the conditions $\partial_i P=0$,
$i=1,2,\cdots ,p$, follow from the construction of polynomials
of the space $\tilde{\cal H}^{({\bf t})}_{s,s'}$. In order to prove
the fulfillment of the conditions $\bar\partial_i P=0$, $i=1,2,\cdots ,p$,
we note that according to formulas (8) and (9) the space
$\tilde{\cal H}^{({\bf t})}_{s,s'}$ is elementwise invariant with respect
to the subalgebra $U_q({\rm gl}_{p})$. Moreover, this space is invariant
and irreducible with respect to the subalgebra $U_q({\rm gl}_{n-p})$
acting on ${\bf t}$.

Now we rearrange elements $z_{p+1},\cdots ,z_{n}, w_{p+1},\cdots ,w_{n}$
in each of polynomials of $\tilde{\cal H}^{({\bf t})}_{s,s'}$ such that
in each summand (monomial) elements $w_{p+1},\cdots ,w_{n}$ stand before the
elements $z_{p+1},\cdots ,z_{n}$. We denote the space
$\tilde{\cal H}^{({\bf t})}_{s,s'}$ with this rearrangement in polynomials
by ${\cal L}^{({\bf t})}_{s,s'}$. Because of elementwise invariance
with respect to $U_q({\rm gl}_{p})$, the space ${\cal L}^{({\bf t})}_{s,s'}$
can be represented as a direct sum
$$
{\cal L}^{({\bf t})}_{s,s'}={\cal R}_{s,s'}\oplus Q_{\bf y}{\cal
R}_{s-1,s'-1}
\oplus Q_{\bf y}^2 {\cal R}_{s-2,s'-2}\oplus \cdots , \eqno (51)
$$
where ${\cal R}_{s-j,s'-j}$ denote the space of homogeneous polynomials
in which $w_{p+1},\cdots ,w_{n}$ stand before $z_{p+1},\cdots ,z_{n}$.
Due to formulas (8) and (9), the spaces ${\cal R}_{s-j,s'-j}$ are invariant
with respect to $U_q({\rm gl}_{n-p})$. However, the representation of
$U_q({\rm gl}_{n-p})$ on ${\cal L}^{({\bf t})}_{s,s'}$ is irreducible.
Therefore, the decomposition (51) contains only one summand and
${\cal L}^{({\bf t})}_{s,s'}={\cal R}_{s,s'}$. It is clear that for
elements of ${\cal R}_{s,s'}$ the conditions $\bar\partial_i P=0$,
$i=1,2,\cdots ,p$, are fulfilled. Lemma is proved.

{\it Corollary 1:} {\it Elements $P$ of the space
$\tilde{\cal H}^{({\bf t})}_{s,s'}$ satisfy the relation
$\Delta _{({\bf t})}P=0$.}

{\it Corollary 2:} {\it Elements $P$ of the space
$\tilde{\cal H}^{({\bf t})}_{s,s'}$ are $q$-harmonic, that is,
$\Delta _q P=0$.}

Corollary 1 follows from (47)--(49). Corollary 2 follows from Corollary
1 and formula (46).

{\it Lemma 3:} {\it For polynomial $h_{s,s'}({\bf t})\in
\tilde{\cal H}^{({\bf t})}_{s,s'}$
and arbitrary polynomial $f({\bf y})$ we have}
$$
\hat\Delta _{({\bf y})} h_{s,s'}({\bf t})f({\bf y})=q^{s-s'}h_{s,s'}({\bf
t})
\hat\Delta _{({\bf y})}f({\bf y}).
$$

{\it Proof:} We first prove the relations
$\partial_i h_{s,s'}({\bf t})f({\bf y})=q^{s'}h_{s,s'}({\bf t})\partial_i
f({\bf y})$, $i=1,\cdots,p$.
The polynomial $h_{s,s'}({\bf t})$ can be represented in the form of a
linear combination of monomials $z_{p+1}^{k_{p+1}}\cdots z_{n}^{k_n}
w_{n}^{l_n}\cdots w_{p+1}^{l_{p+1}}$,
where $k_{p+1}+\cdots+k_n=s$, $l_{p+1}+\cdots+l_n=s'$.
We have
$$
\partial_i z_{p+1}^{k_{p+1}}\cdots z_{n}^{k_n}w_{n}^{l_n}\cdots
w_{p+1}^{l_{p+1}} f({\bf y})=
z_{p+1}^{k_{p+1}}\cdots z_{n}^{k_n} (\partial_i w_{n}^{l_n}\cdots
w_{p+1}^{l_{p+1}} f({\bf y}))
$$
$$
=q^{s'} z_{p+1}^{k_{p+1}}\cdots z_{n}^{k_n}w_{n}^{l_n}\cdots
w_{p+1}^{l_{p+1}} (\partial_i f({\bf y})),
$$
where the relation $\partial_j f({\bf y})=0$ and relations from
section III were used.
It proves our relations. We analogously prove the relations
$\bar\partial_i h_{s,s'}({\bf t})f({\bf y})=q^{-s}h_{s,s'}({\bf t})
\bar\partial_i f({\bf y})$,
$i=1,\cdots,p$. In this case, it is useful to represent the polynomial
$h_{s,s'}({\bf t})$ in the form of a linear combination
of monomials $w_{p+1}^{l_{p+1}}\cdots w_{n}^{l_n}z_{n}^{k_n}\cdots
z_{p+1}^{k_{p+1}}$
(such representation is possible due to Lemma 2). Now the lemma follows
from explicit formula for
$\hat\Delta _{({\bf y})}$. Lemma is proved.

Since $\Delta_{({\bf t})} \left( Q_{\bf y}^u h_{s,s'}({\bf t}) h_{r,r'}
({\bf y}) \right)=0$,
then using Lemma 3 and relation (18) with $n$ replaced by $p$ we have
$$
\Delta_q \left( Q_{\bf y}^u h_{s,s'}({\bf t}) h_{r,r'}({\bf y}) \right)=
\hat\Delta_{({\bf y})}
\left( Q_{\bf y}^u h_{s,s'}({\bf t}) h_{r,r'}({\bf y}) \right)=
q^a h_{s,s'}({\bf t}) \Delta_{({\bf y})} Q_{\bf y}^u h_{r,r'}({\bf y})
$$
$$=q^a[u][p+u+r+r'-1]\left( Q_{\bf y}^{u-1} h_{s,s'}({\bf t})
h_{r,r'}({\bf y}) \right),
$$
where $a=2(s-s')u+s'-s$.

Now we may find a harmonic projection of the polynomials (50).
Denoting this projection
by $h_{m,m'}^{(r,r';s,s')}({\bf z},{\bf w})$ we have
$$
h_{m,m'}^{(r,r';s,s')}({\bf z},{\bf w})=\sum _{k=0}^{{\rm min}\
(m,m')} \alpha _k Q^k \Delta_q ^k \left( Q_{\bf y}^u h_{s,s'}({\bf
t}) h_{r,r'}({\bf y}) \right)
$$      $$
=\left( \sum _{k=0}^{u}
\alpha _k Q^k q^{(s-s'+p-1)k}\frac{[u]![r+r'+p+u-1]!}{[u-k]!
[r+r'+p+u-k-1]!}
Q_{\bf y}^{u-k} \right) h_{s,s'}({\bf t}) h_{r,r'}({\bf y}) ,
$$
where  $\alpha _k$ is determined by formula
(33). Denoting the expression in the parentheses by
$t^{n,p;m,m'}_{r,r';s,s'}(Q_{\bf y},Q_{\bf t})$, we have
$$
h_{m,m'}^{(r,r';s,s')}({\bf z},{\bf w})={\sf H}_{m,m'}\left( Q_{\bf
y}^{m-r-s} h_{s,s'}({\bf t}) h_{r,r'}({\bf y})\right) =
t^{n,p;m,m'}_{r,r';s,s'}(Q_{\bf y},Q_{\bf t}) h_{s,s'}({\bf t})
h_{r,r'}({\bf y}) . \eqno (52)
$$

After some simple transformations, we obtain for
$t^{n,p;m,m'}_{r,r';s,s'}(Q_{\bf y},Q_{\bf t})$ the expression
$$
t^{n,p;m,m'}_{r,r';s,s'}(Q_{\bf y},Q_{\bf t})=Q_{\bf y}^u \sum
_{k=0}^u \frac{(q^{-2u};q^2)_k
(q^{-2(r+r'+p+u-1)};q^2)_k}{(q^{-2(m+m'+n-2)};q^2)_k (q^2;q^2)_k}
q^{k\sigma}Q^kQ_{\bf y}^{-k},
$$
where $\sigma =-2n-2s'+2+2p$. Taking into account the
definition of the basis hypergeometric function ${}_2 \varphi _1$, we
derive
$$
t^{n,p;m,m'}_{r,r';s,s'}(Q_{\bf y},Q_{\bf t})= Q_{\bf y}^u\, {}_2
\varphi _1(q^{-2u},q^{-2(r+r'+p+u-1)};\ q^{-2(m+m'+n-2)};\ q^2,\
QQ^{-1}_{\bf y} q^\sigma ).
$$
Applying the relation
$$
{}_2 \varphi _1 (q^{-n},b;\ c;\ q,z)=q^{-(n+1)n/2}(-z)^n
\frac{(b;q)_n}{(c;q)_n} {}_2 \varphi _1 (q^{-n}, q^{1-n}/c;\
q^{1-n}/b;\ q, cq^{n+1}/bz)
$$
(see, for example, formula (2) of Section 14.1.8 in Ref. 14)
we reduce this expression to
$$
t^{n,p;m,m'}_{r,r';s,s'}(Q_{\bf y},Q_{\bf t})=(-q^\sigma )^uq^{-(u+1)u}
\frac{(q^{-2(r+r'+p+u-1)};q^2)_u}{(q^{-2(m+m'+n-2)};q^2)_u}
Q^u
$$        $$
\times {}_2 \varphi _1(q^{-2u},q^{2(m+m'+n-u-1)};\
q^{2(r+r'+p)};\ q^2,\ q^{-2s+2} Q_{\bf y}/Q ).
$$
Using the definition of the little $q$-Jacobi polynomials, we
derive from here that
$$
t^{n,p;m,m'}_{r,r';s,s'}(Q_{\bf y},Q_{\bf t})=(-q^\sigma )^uq^{-(u+1)u}
\frac{(q^{-2(r+r'+p+u-1)};q^2)_u}{(q^{-2(m+m'+n-2)};q^2)_u}
Q^u
$$        $$
\times P_u^{(r+r'+p-1,s+s'+n-p-1)}( q^{-2s} Q_{\bf y}/Q;\
q^2 ). \eqno (53)
$$
Thus, we proved that {\it the projection ${\sf H}_{m,m'}\left( Q_{\bf
y}^{m-r-s} h_{s,s'}({\bf t}) h_{r,r'}({\bf y})\right)$ is given by
formula (52), where $t^{n,p;m,m'}_{r,r';s,s'}$ is determined by
(53).} The restriction $\tau h^{(r,r';s,s')}_{m,m'} ({\bf
z},{\bf w})$ of this projection onto the quantum sphere
$S^{\Bbb C}_{q,n-1}$ is given by
$$
\tau h^{(r,r';s,s')}_{m,m'}({\bf z},{\bf w})=(
\tau t^{n,p;m,m'}_{r,r';s,s'})(Q_{\bf y})h_{s,s'}({\bf t})
h_{r,r'}({\bf y}) ,
$$
where
$(\tau t^{n,p;m,m'}_{r,r';s,s'})(Q_{\bf y})=c
P_u^{(r+r'+p-1,s+s'+n-p-1)}( q^{-2s} Q_{\bf y}/Q; q^2 )$
($c$ is the multiplier from the right hand side of (53)).

For the scalar product of polynomials of the form (52) we have
$$
\langle h^{(r,r';s,s')}_{m,m'},h^{({r''},{r'''};
{s''},{s'''})}_{m,m'} \rangle =0\ \ \ \ {\rm if}\ \ \ \
(r,r',s,s')\ne (r'',r''',s'',s''')
$$
(since the spaces ${\cal H}^{({\bf y})}_{r,r'}$ and
${\cal H}^{({\bf y})}_{{r''},{r}'''}$ and the spaces
$\tilde{\cal H}^{({\bf t})}_{s,s'}$ and
$\tilde{\cal H}^{({\bf t})}_{{s''},{s'''}}$ are orthogonal). If
$(r,r',s,s')= (r'',r''',s'',s''')$, then
the norm of the polynomial (52) reduces to the orthogonality relation
for $q$-Jacobi polynomials and to norms of
$h_{s,s'}({\bf t})$ and $h_{r,r'}({\bf y})$.

In order to obtain a $q$-analogue of separation of variables in this
case we have to take bases of the spaces ${\cal H}^{({\bf y})}_{r,r'}$
and $\tilde{\cal H}^{({\bf t})}_{s,s'}$ in separated coordinates (as it was
made in section IX).
\bigskip

\noindent{\sf ACKNOWLEDGMENT}
\medskip

The research by the first author (NZI) was partially supported by INTAS
grant
No.~2000-334.

\bigskip

\noindent ${}^1$N. Ja. Vilenkin and A. U. Klimyk, {\it
Representation of Lie Groups and Special Functions} (Kluwer
Academic, Dordrecht, 1993), Vol. 2.

\noindent ${}^2$W. Miller, {\it Symmetry and Separation of
Variables} (Addison-Wesley, Massachusetts, 1978).

\noindent ${}^3$P. G. A. Floris, Compositio Math. {\bf 108}, 123
(1997).

\noindent ${}^4$P. G. A. Floris, {\it On quantum groups,
hypergroups and $q$-special functions} (Ph. D. Theses,
Leiden University, 1995).

\noindent ${}^5$N. Z. Iorgov and A. U. Klimyk, J. Math. Phys. {\bf
42}, 1326 (2001).

\noindent ${}^6$A. Klimyk and K. Schm\"udgen, {\it Quantum Groups
and Their Representations} (Springer, Berlin, 1997).

\noindent ${}^7$N. Ya. Reshetikhin, L. A. Takhtajan, and L. D.
Faddeev, Leningrad Math. J. {\bf 1}, 193 (1990).

\noindent ${}^8$M. Noumi, H. Yamada, and K. Mimachi, Japan J.
Math. {\bf 19}, 31 (1993).

\noindent ${}^9$Ch.-S. Chu, P.-M. Ho, and B. Zumino,
arXiv:q-alg/9510021.

\noindent ${}^{10}$B. Kostant, Am. J. Math. {\bf 83}, 327 (1963).

\noindent ${}^{11}$I. M. Burban and A. U. Klimyk, J. Phys. A {\bf
26}, 2139 (1993).

\noindent ${}^{12}$G. Gasper and M. Rahman, {\it Basic
Hypergeometric Functions} (Cambridge University Press, Cambridge,
1990).

\noindent ${}^{13}$G. E. Andrews, R. Askey, and R. Roy, {\it
Special Functions} (Cambridge University Press, Cambridge, 1999)

\noindent ${}^{14}$N. Ja. Vilenkin and A. U. Klimyk, {\it
Representation of Lie Groups and Special Functions} (Kluwer
Academic, Dordrecht, 1992), Vol. 3.

\end{document}